\documentclass[letter,10pt]{amsart}
\usepackage{amsmath,amssymb,colordvi,verbatim}
\usepackage{multicol}
\usepackage[utf8]{inputenc}
\usepackage[pagewise]{lineno}
\newtheorem{theorem}{Theorem}

\newtheorem{lemma}[theorem]{Lemma}

\newtheorem{proposition}[theorem]{Proposition}
\newtheorem{definition}[theorem]{Definition}

\newtheorem{remark}[theorem]{\it Remark}

\newcommand{\C}{\mathbb{C}}

\newcommand{\N}{\mathbb{N}}

\makeatletter
\@namedef{subjclassname@2020}{%
  \textup{2020} Mathematics Subject Classification}
\makeatother

\begin{document}

\title[Order in Abstract fractional differential systems]
{Determination of the order in Abstract fractional differential equations.}

\author{Rodrigo Ponce}
\address{Universidad de Talca, Instituto de Matem\'aticas, Casilla 747, Talca-Chile.}
\email{rponce@inst-mat.utalca.cl,rponce@utalca.cl}

\subjclass[2020]{Primary 34K29, Secondary 34A55, 47D06, 26A33}

\date{March 15, 2023 and, in revised form, XXXX, XXXX}

\keywords{Resolvent families; inverse problems; fractional differential equations; unbounded linear operators}


\begin{abstract}
In this paper we identify, for small $t$ and a fixed $T>0,$ the order $\alpha>0$ in the abstract fractional differential equation
$$\partial^\alpha u(t)=Au(t),$$
where the time-fractional derivative $\partial^\alpha$ is understood in the sense of Caputo and Riemann-Liouville, $A$ is a closed (possibly unbounded) linear operator in a Banach space $X,$ and $0<\alpha<1$ or $1<\alpha<2.$
\end{abstract}

\maketitle

\section{Introduction}

The problem of finding or approximating the order in time-fractional differential equations has been widely studied in the last ten years. See for instance \cite{Al-As-20,As-Um-20,Ha-Na-Wa-Ya-13,Ji-Ka-21,Ji-Ru-15,Ka-Ru-21,Li-Hu-Ya-20,Li-Zi-Ya-19,Li-Zh-20,Lu-11,Ya-21,Ze-Ch-Wa-19}. One of the most notable contributions is the paper \cite{Ha-Na-Wa-Ya-13}, where authors consider (for $0<\alpha<1$) the fractional differential equation for the Caputo fractional derivative
\begin{equation}\label{EqIntro1}
\partial^\alpha_t u(x,t)=Au(x,t), \quad x\in\Omega, t>0
\end{equation}
under the initial condition $u(x,0)=u_0(x), x\in \Omega,$ where $\Omega$ and $A$ are defined as follows: For a bounded open set  $\Omega\subset\mathbb{R}^N$ with sufficiently smooth boundary $\partial \Omega,$ let $X$ be the Hilbert space $L^2(\Omega).$ On $X,$ the operator $\mathcal{A}$ is defined by $\mathcal{A} u(x)=\sum_{i=1}^N \frac{\partial }{\partial x_i}\left(\sum_{j=1}^N A_{ij}(x)\frac{\partial }{\partial x_i}u(x)\right),  u\in X,$ where $A_{ij}=A_{ji}$ for any $1\leq i,j\leq N.$ Suppose that there exists a constant $\gamma>0$ such that $\sum_{i,j=1}^N A_{ij}(x)\xi_i\xi_j\geq \gamma|\xi|^2,$
for all $\xi\in\mathbb{R}^N$ and $x\in \overline{\Omega}.$ The operator $A:D(A)\to X$ is defined by
\begin{equation*}
(Au)(x)=(\mathcal{A}u)(x),\quad x\in \Omega,
\end{equation*}
where $D(-A)=H^2(\Omega)\cap H_0^1(\Omega).$ The operator $-A$ has a discrete spectrum and its eigenvalues satisfy its eigenvalues satisfy $0<\lambda_1\leq\lambda_2\leq\cdots\le\lambda_n\le\cdots$ and $\lim_{n\to\infty}\lambda_n=\infty.$ Now, if  $\phi_n\in H^2(\Omega)\cap H_0^1(\Omega)$ denotes the normalized eigenfunction associated with $-\lambda_n,$ then, by the Fourier method (see \cite{Sa-Ya-11}), the solution $u$ to \eqref{EqIntro1} is given by
\begin{equation}\label{EqIntro2}
u(x,t)=\sum_{n=1}^\infty\langle u_0,\phi_n\rangle_{L^2(\Omega)} E_{\alpha,1}(-\lambda_n t^\alpha)\phi_n(x),\quad x\in \Omega, t>0,
\end{equation}
where for any $\alpha,\beta>0$ and $z\in\C,$ $E_{\alpha,\beta}$ denotes the Mittag-Leffler function which is defined by
$ E_{\alpha,\beta}(z):=\sum_{k=0}^\infty\frac{z^k}{\Gamma(\alpha k+\beta)}.$
If $u(t,x)$ denotes the solution to \eqref{EqIntro1}, $u_0\in C_0^\infty (\Omega)$ with $Au_0(x)\neq 0$ and $x_0$ is a fixed element in $\Omega$, then the order $\alpha$ in \eqref{EqIntro1} is given by (see \cite[Theorem 1]{Ha-Na-Wa-Ya-13})
\begin{equation}\label{EqIntro5}
\alpha =\lim_{t\to 0^+}\frac{t\frac{\partial u}{\partial t}(x_0,t)}{u(x_0,t)-u_0(x_0)}.
\end{equation}
A similar result holds for $t\to +\infty$ (see also \cite[Theorem 1]{Ha-Na-Wa-Ya-13}). Thus, to determinate the order $\alpha$ we need to know $u(x_0,t)$ and $\frac{\partial u}{\partial t}(x_0,t)$ for $t>0$ on an interval close to $0$ (or $+\infty$). As the authors mention in \cite[p. 440]{Li-Li-Ya-19}, "the problems of the recovery of the fractional orders are far from satisfactory since all the publications either assumed the homogeneous boundary condition or studied this inverse problem by the measurement in $t\in (0,\infty)."$ Therefore, the problem of finding the order $\alpha$ in \eqref{EqIntro1} in terms of its solution $u(x,t)$ for a fixed time $t>0$ remains as an open problem.



Now, if for any $t\geq 0,$ we define the family of linear operators $S_{\alpha,\beta}(t):X\to X$  by
\begin{equation*}
S_{\alpha,\beta}(t)u(x):=\sum_{n=1}^\infty \langle u_0,\phi_n\rangle_{L^2(\Omega)} t^{\beta-1}E_{\alpha,\beta}(-\lambda_n t^\alpha)\phi_n(x), \quad x\in \Omega, u\in X,
\end{equation*}
then, the solution \eqref{EqIntro2} to equation \eqref{EqIntro1} can be written as
\begin{equation}\label{EqIntro3}
u(x,t)=S_{\alpha,1}(t)u_0(x),\quad x\in\Omega, t\geq 0.
\end{equation}
The properties of the Laplace transform of the Mittag-Leffler function imply that $\{S_{\alpha,\beta}(t)\}_{t\geq 0}$ corresponds to an $(\alpha,\beta)$-fractional resolvent family generated by $A,$ see for instance \cite{Li-19}. This theory allows us write the solutions to fractional differential equations (for the Caputo and Riemann-Liouville derivatives) in case $0<\alpha<1$ and $1<\alpha<2$ as a variation of parameters formulas. In fact, consider the fractional differential equations for the Caputo fractional derivative,
\begin{equation}\label{eqCaputoI}
\left\{\begin{array}{lcl}
\partial_t^\alpha u(x,t) &=& Au(x,t), \quad t>0\\
u(x,t)&=&0,\qquad \quad x\in \partial\Omega, t>0\\
u(x,0)&=&u_0(x),\quad x\in\Omega,\\
\end{array}\right.
\end{equation}
(for $0<\alpha<1$) and
\begin{equation}\label{eqCaputoII}
\left\{\begin{array}{lcl}
\partial_t^\alpha u(x,t) &=& Au(x,t), \quad t>0\\
u(x,t)&=&0,\qquad \quad x\in \partial\Omega, t>0\\
u(x,0)&=&u_0(x),\quad x\in\Omega\\
\partial_tu(x,0)&=&u_1(x),\quad x\in\Omega,\\
\end{array}\right.
\end{equation}
(for $1<\alpha<2$) where $u_0,u_1\in X.$ By \cite{Sa-Ya-11}, the solution to \eqref{eqCaputoI} is given by \eqref{EqIntro3} and the solution to \eqref{eqCaputoII} is
\begin{equation*}
u(x,t)=\sum_{n=1}^\infty\langle u_0,\phi_n\rangle_{L^2(\Omega)} E_{\alpha,1}(-\lambda_n t^\alpha)\phi_n(x)+\sum_{n=1}^\infty\langle u_1,\phi_n\rangle_{L^2(\Omega)} tE_{\alpha,2}(-\lambda_n t^\alpha)\phi_n(x),
\end{equation*}
and therefore, the solutions to \eqref{eqCaputoI} and \eqref{eqCaputoII} can be written, in terms of the resolvent family $\{S_{\alpha,\beta}(t)\}_{t\geq 0},$ respectively, as
\begin{equation*}
u(x,t)=S_{\alpha,1}(t)u_0(x),\quad u(x,t)=S_{\alpha,1}(t)u_0(x)+S_{\alpha,2}(t)u_1(x),\quad x\in \Omega, t\geq 0.
\end{equation*}

Now, if we consider the fractional differential equations for the Riemann-Liouville fractional derivatives
\begin{equation}\label{eqRLI}
\left\{\begin{array}{lcl}
^R\partial_t^\alpha u(x,t) &=& Au(x,t), \quad t>0\\
u(x,t)&=&0,\qquad \quad x\in \partial \Omega, t>0\\
(g_{1-\alpha}\ast u)(x,0)&=&u_0(x),\quad x\in\Omega,\\
\end{array}\right.
\end{equation}
(for $0<\alpha<1$) and
\begin{equation}\label{eqRLII}
\left\{\begin{array}{lcl}
^R\partial_t^\alpha u(x,t) &=& Au(x,t), \quad t>0\\
u(x,t)&=&0,\qquad \quad x\in \partial \Omega, t>0\\
(g_{2-\alpha}\ast u)(x,0)&=&u_0(x),\quad x\in\Omega\\
\partial_t(g_{2-\alpha}\ast u)(x,0)&=&u_1(x),\quad x\in\Omega,\\
\end{array}\right.
\end{equation}
(for $1<\alpha<2$), where $u_0,u_1\in X,$ $^R\partial_t^\alpha$ corresponds to the Riemann-Liouville fractional derivative and $g_{2-\alpha}(t):=t^{1-\alpha}/\Gamma(2-\alpha),$ then, the solutions to \eqref{eqRLI} and \eqref{eqRLII} are given respectively, by (see for instance  \cite{Pon-20}),
\begin{equation*}
u(x,t)=\sum_{n=1}^\infty\langle u_0,\phi_n\rangle_{L^2(\Omega)} t^{\alpha-1}E_{\alpha,\alpha}(-\lambda_n t^\alpha)\phi_n(x),
\end{equation*}
and
\begin{equation*}
u(x,t)=\sum_{n=1}^\infty\langle u_0,\phi_n\rangle_{L^2(\Omega)} t^{\alpha-2}E_{\alpha,\alpha-1}(-\lambda_n t^\alpha)\phi_n(x)+\sum_{n=1}^\infty\langle u_1,\phi_n\rangle_{L^2(\Omega)} t^{\alpha-1}E_{\alpha,\alpha}(-\lambda_n t^\alpha)\phi_n(x),
\end{equation*}
which can be written, respectively, as
\begin{equation*}
u(x,t)=S_{\alpha,\alpha}(t)u_0(x),\quad u(x,t)=S_{\alpha,\alpha-1}(t)u_0(x)+S_{\alpha,\alpha}(t)u_1(x).
\end{equation*}

The resolvent families have been extensively studied, both in abstract settings and in applications (see for instance \cite{Ca-Pl-15,Cu-Pa-04,Ei-Ko-04,He-Li-Zo-20,He-Me-Po-21,Li-19,Pon-20,Wa-Ch-Xi-12}). The operators $\{S_{\alpha,\beta}(t)\}_{t\geq 0}$ are well-known in some cases: the uniqueness of the Laplace transform implies that $\{S_{1,1}(t)\}_{t\geq 0}$ is the  $C_0$-semigroup generated by $A,$ $\{S_{2,1}(t)\}_{t\geq 0}$ and $\{S_{2,2}(t)\}_{t\geq 0}$ are, respectively, the cosine and sine family generated by $A,$ see \cite{Ar-Ba-Hi-Ne-11}. Now, for $1\leq \alpha\leq 2$ and $\beta=1,$ $\{S_{\alpha,1}(t)\}_{t\geq 0}$ is an $\alpha$-times resolvent \cite{Li-Ch-Li-10}, the case $1\leq \alpha=\beta\leq 2$ corresponds to an $\alpha$-order resolvent (see \cite{Li-Pe-Ji-12}) and if $\alpha=1$ and $\beta=n+1, n\in\mathbb{N},$ then we get an $n$-times integrated semigroup, see \cite{Ar-Ba-Hi-Ne-11}.

In this paper, we explore the resolvent families $\{S_{\alpha,\beta}(t)\}_{t\geq 0}$ to identify the order $\alpha$ (for small times $t>0$ and a fixed time $T>0$) in the fractional differential equations \eqref{eqCaputoI}--\eqref{eqRLII}, where $A$ is a closed linear operator in a Banach space $X.$ More specifically, we consider the abstract fractional differential equation for the Caputo fractional derivative
\begin{equation}\label{EqIntro4}
\left\{\begin{array}{lcl}
\partial_t^\alpha v(t) &=& Av(t), \quad t>0\\
v(0)&=&x,\\
\end{array}\right.
\end{equation}
where $0<\alpha<1,$ $A$ is a closed linear operator defined in a Banach space $X$ and $x\in X.$ If $v(t;x)=S_{\alpha,1}(t)x$ is the solution to Problem \eqref{EqIntro4} (where $\{S_{\alpha,1}(t)\}_{t\geq 0}$ is the $(\alpha,1)$-resolvent family generated by $A$), then we prove in Theorem \ref{Th3.1} that
\begin{equation*}
\alpha x=\lim_{t\to 0^+} t\psi_t\varphi^{-1}_t(x),
\end{equation*}
where  $\varphi_t:X\to X$ and $\psi_t:X\to X,$ are defined respectively, by $\varphi_t(x)=v(t;x)-x$ and $\psi_t(x)=v'(t;x).$ This is exactly the abstract version of formula \eqref{EqIntro5}. Moreover, if $T>0$ is a fixed time, $A$ generates the $(\alpha,1)$-resolvent family $\{S_{\alpha,1}(t)\}_{t\geq 0}$ and $x$ is an element in the Banach space $X,$ then the order $\alpha$ verifies
\begin{equation*}
Tv(T;x)-(g_1\ast v)(T;x)=\alpha \left[(S_{\alpha,1}\ast v)(T;x)-(g_1\ast v)(T;x)\right],
\end{equation*}
where $v(t;x)$ is the solution to Problem \eqref{EqIntro4}, see Theorem \ref{Th5.1} and Remark \ref{Rem1}. Thus, to determinate $\alpha$ we only to know the solution $v(T;x),$ its integral $(g_1\ast v)(T;x)$ and the convolution $(S_{\alpha,1}\ast v)(T;x)$ for any fixed $x\in X$ and $T>0.$ This implies that
in Equation \eqref{EqIntro1}, the order $\alpha$ verifies
\begin{equation*}
\alpha=\frac{Tu(x_0,T)-(g_1\ast u)(x_0,T)}{(S_{\alpha,1}\ast u)(x_0,T)-(g_1\ast u)(x_0,T)}
\end{equation*}
for any fixed $x_0\in \Omega$ and $T>0$ such that $(S_{\alpha,1}\ast u)(x_0,T)-(g_1\ast u)(x_0,T)\neq 0,$ where $u(x,t)$ is given by \eqref{EqIntro3}. This result gives an answer to the problem proposed in \cite[p. 440]{Li-Li-Ya-19} and allow us to find the order  $\alpha$ in Equation \eqref{EqIntro1} in terms of the solution $u(x,t)$ for a fixed time $t>0.$

Moreover, we obtain here similar results for \eqref{EqIntro4} in case $1<\alpha<2$ and for the abstract fractional differential equation for the Riemann-Liouville fractional derivative in case $0<\alpha<1$ and $1<\alpha<2.$

The paper is organized as follows. In Section \ref{Sect2} we give the preliminaries on fractional calculus and fractional resolvent families generated by a closed linear operator $A.$ In Section \ref{Sect3} we identify $\alpha\in (0,1)$ for small times for the Caputo and Riemann-Liouville fractional derivatives. Section \ref{Sect4} is devoted to the same problem, but with $\alpha\in (1,2).$ In Sections \ref{Sect5} and \ref{Sect6} we identify, respectively, $\alpha\in (0,1)$ and $\alpha\in (1,2),$ for a fixed time $T>0$ for the Caputo and Riemann-Liouville fractional derivatives. Finally, we illustrate our results with some examples.

\section{Preliminaries}\label{Sect2}

Let $X\equiv(X,\|\cdot\|)$ be a Banach space. By $\mathcal{B}(X)$ we denote the Banach space of all bounded and linear operators from $X$ into $X.$ For a given closed linear operator $A$ on $X,$ $\rho(A)$ denotes its resolvent set and  $R(\lambda,A)=(\lambda-A)^{-1}$ its resolvent operator, which is defined for all $\lambda\in\rho(A).$

A strongly continuous family of linear operators $\{S(t)\}_{t\geq 0}\subset \mathcal{B}(X)$ is called {\em exponentially bounded} if there exist constants $M>0$ and $w\in\mathbb R$ such that $\|S(t)\|\leq Me^{wt},$ for all $t>0.$

For a given $\alpha>0,$ we define the function $g_\alpha$ as $g_\alpha(t):=\frac{t^{\alpha-1}}{\Gamma(\alpha)},$ where $\Gamma(\cdot)$ stands the Gamma function. It is easy to see that if $\alpha,\beta>0,$ then the functions $g_\gamma$ satisfy the semigroup law $g_{\alpha+\beta}(t)=(g_\alpha\ast g_\beta)(t),$ where $(f\ast g)$ denotes the finite convolution $(f\ast g)(t):=\int_0^t f(t-s)g(s)ds.$

For $n-1<\alpha<n,$ where $n\in\N,$ the Caputo and Riemann-Liouville fractional derivatives of order $\alpha$ of a function $f$ are defined, respectively, by
\begin{equation*}
\partial_t^\alpha f(t):=(g_{n-\alpha}\ast f')(t)=\int_0^t g_{n-\alpha}(t-s)f^{(n)}(s)ds, \quad ^{R}\partial_t^\alpha f(t):=\frac{d^n}{dt^n}\int_0^t g_{n-\alpha}(t-s)f(s)ds.
\end{equation*}
If $\alpha=1$ or $\alpha=2,$ then $\partial^1_t= \,^R\partial_t^1=\frac{d}{dt}$ and $\partial^2_t= \,^R\partial_t^2=\frac{d^2}{dt^2}.$  For more details, examples and applications on fractional calculus, we refer to the reader to \cite{Ko-Lu-19}. In this paper, our focus is in $0<\alpha<1$ and $1<\alpha<2.$

\begin{definition}\label{def3.1}
Let $A$ be a closed and linear operator defined on a Banach space $X.$ Given $\alpha,\beta>0$ we say that $A$ is the generator of an $(\alpha,\beta)$-resolvent family, if there exist $\omega\geq 0$ and a strongly continuous function $S_{\alpha,\beta}: (0,\infty)\to \mathcal{B}(X)$ such that $S_{\alpha,\beta}(t)$ is exponentially bounded, $\left\{\lambda^\alpha: {\rm Re} \lambda>\omega\right\}\subset \rho(A),$ and for all $x\in X,$
\begin{equation}\label{eqRelovent-1}
\lambda^{\alpha-\beta}\left(\lambda^\alpha-A\right)^{-1}x=\int_0^\infty e^{-\lambda t}S_{\alpha,\beta}(t)xdt, \,\,\, {{\rm Re} \lambda>\omega}.
\end{equation}
In this case, $\{S_{\alpha,\beta}(t)\}_{t>0}$ is called the {\em $(\alpha,\beta)$-resolvent family} generated by $A.$
\end{definition}

Moreover, if an operator $A$ with domain $D(A)$ is the infinitesimal generator of $S_{\alpha,\beta}(t),$ then, for $x\in D(A)$, we have
\begin{eqnarray*}
Ax=\lim_{t\to 0^+}\frac{S_{\alpha,\beta}(t)x-g_\beta(t)x}{g_{\alpha+\beta}(t)}.
\end{eqnarray*}

We notice that the case $S_{1,1}(t)$ corresponds to a $C_0$-semigroup, $S_{2,1}(t)$ is a cosine family and $S_{2,2}(t)$ is a sine family generated by $A.$ In the scalar case, that is, when $A=\rho I,$ where $\rho\in \mathbb{C}$ and $I$ denotes the identity operator, we have, by the uniqueness of the Laplace transform, that $S_{\alpha,\beta}(t)$ corresponds to the function $t^{\beta-1}E_{\alpha,\beta}(\rho t^\alpha).$
Finally, for $0<\alpha<1$ and $\beta\geq \alpha,$ let $\{S_{\alpha,\beta}(t)\}_{t\geq 0}$ be the family of operators defined by
\begin{eqnarray*}
S_{\alpha,\beta}(t)f(s):=\int_{0}^s f(s-r)\varphi_{\alpha,\beta-\alpha}(t,r)dr,
\end{eqnarray*}
where $s\in \mathbb{R}_+,$ $f\in L^1(\mathbb{R}_+)$ and the function $\varphi_{a,b}(t,r)$ is defined by
\begin{eqnarray*}\label{eqLevyPr}
\varphi_{a,b}(t,r):=t^{b-1}W_{-a,b}(-rt^{-a}),\quad a>0, b\geq 0,
\end{eqnarray*}
where $W_{-a,b}(z):=\sum_{n=0}^\infty \frac{z^n}{n!\Gamma(-an+b)}$ ($z\in\mathbb{C}$) denotes the Wright function. Then, $\{S_{\alpha,\beta}(t)\}_{t\geq 0}$ is an $(\alpha,\beta)$-resolvent family on the Banach space $X=L^1(\mathbb{R}_+)$ generated by $A=-\frac{d}{dt}.$ See \cite[Example 11]{Ab-Mi-15}.

From \cite{Ab-Al-20,Ab-Mi-15} or \cite{Li-00} we have  the following result that gives some important properties of the resolvent family $\{S_{\alpha,\beta}(t)\}_{t\geq 0}.$
\begin{proposition}\label{Proposition1}
If $\alpha,\beta>0$ and  $A$ generates an $(\alpha,\beta)$-resolvent family $\{S_{\alpha, \beta}(t)\}_{t>0},$ then
\begin{enumerate}
  \item $\displaystyle\lim_{t\to 0^+}\frac{S_{\alpha, \beta}(t)x}{g_\beta(t)}=x,$ for all $x\in X.$
  \item $S_{\alpha, \beta}(t)x\in D(A)$ and $S_{\alpha, \beta}(t)Ax=AS_{\alpha, \beta}(t)x$ for all $x\in D(A)$ and $t>0.$
  \item For all $x\in D(A),$
  \begin{equation*}
  S_{\alpha,\beta}(t)x=g_\beta(t)x+\int_0^t g_\alpha(t-s)AS_{\alpha, \beta}(s)xds.
  \end{equation*}
  \item $\int_0^t g_\alpha(t-s)S_{\alpha, \beta}(s)xds\in D(A)$ and
        \begin{equation*}
         S_{\alpha,\beta}(t)x=g_\beta(t)x+A\int_0^t g_\alpha(t-s)S_{\alpha, \beta}(s)xds,
        \end{equation*}
      for all $x\in X.$
\end{enumerate}
\end{proposition}

For a given locally integrable function $f:[0,\infty)\to X,$ we define the {\em Laplace transform} of $f,$  denoted by $\hat{f}(\lambda)$ (or $\mathcal{L}(f)(\lambda)$) as
\begin{equation*}
\hat{f}(\lambda)=\int_0^{\infty}e^{-\lambda t}f(t)dt,
\end{equation*}
provided the integral converges for some $\lambda\in\mathbb{C}.$

The following lemmata will be useful for our purposes.

\begin{lemma}\label{Lemma1}
Assume that $A$ is the generator of the family $\{S_{\alpha,\beta}(t)\}_{t\geq 0}.$ If $\gamma>0,$ then $A$ generates the family $\{S_{\alpha,\beta+\gamma}(t)\}_{t\geq 0}$ given by
\begin{equation}
S_{\alpha,\beta+\gamma}(t)=(g_\gamma\ast S_{\alpha,\beta})(t).
\end{equation}
\end{lemma}
\begin{proof}
In fact, for any $\gamma>0$ and $\lambda^\alpha\in\rho(A)$ with ${\rm Re}\lambda>w,$ we have
\begin{equation*}
\mathcal{L}(g_\gamma\ast S_{\alpha,\beta})(\lambda)=\lambda^{\alpha-\beta-\gamma}(\lambda^\alpha-A)^{-1}=\lambda^{\alpha-(\beta+\gamma)}(\lambda^\alpha-A)^{-1}=\hat{S}_{\alpha,\beta+\gamma}(\lambda).
\end{equation*}
And the result follows from the uniqueness of the Laplace transform.
\end{proof}

\begin{lemma}\label{Lemma2}
Assume that $A$ is the generator of the family $\{S_{\alpha,1}(t)\}_{t\geq 0},$ where $\alpha>0.$ Then, the following assertions hold for any $t\geq 0$ and $x\in X,$
\begin{enumerate}
  \item $tS_{\alpha,1}(t)x=-(\alpha-1)(g_1\ast S_{\alpha,1})(t)x+\alpha (S_{\alpha,1}\ast S_{\alpha,1})(t)x.$
  \item $tS_{\alpha,1}'(t)x=\alpha(S_{\alpha,1}'\ast S_{\alpha,1})(t)x.$
\end{enumerate}
\end{lemma}
\begin{proof}
To prove (1), let $h(t):=tS_{\alpha,1}(t)x.$ By the properties of the Laplace transform, we have for any $\lambda^\alpha\in\rho(A),$
\begin{equation*}
\hat{h}(\lambda)=-\frac{d}{d\lambda}(\hat{S}_{\alpha,1}(\lambda))x=-\frac{d}{d\lambda}(\lambda^{\alpha-1}(\lambda^\alpha-A)^{-1}x)=-\frac{(\alpha-1)}{\lambda}\hat{S}_{\alpha,1}(\lambda)x+\alpha\hat{S}_{\alpha,1}(\lambda)\hat{S}_{\alpha,1}(\lambda)x,
\end{equation*}
and the assertion follows from the uniqueness of the Laplace transform. As $S_{\alpha,1}(0)=I,$ to prove the second assertion we only need to note that (2) corresponds to the derivative of (1).
\end{proof}

\begin{lemma}\label{Lemma3}
Assume that $A$ is the generator of the family $\{S_{\alpha,\alpha}(t)\}_{t\geq 0},$ where $0<\alpha<1.$ Then, for any $t\geq 0$ and $x\in X,$
\begin{enumerate}
  \item $tS_{\alpha,\alpha}(t)x=\alpha(g_{1-\alpha}\ast S_{\alpha,\alpha}\ast S_{\alpha,\alpha})(t)x=\alpha(S_{\alpha,1}\ast S_{\alpha,\alpha})(t)x,$ where $S_{\alpha,1}(t)=(g_{1-\alpha}\ast S_{\alpha,\alpha})(t).$
  \item $\int_0^t rS_{\alpha,\alpha}(r)xdr=\alpha[A(g_2\ast S_{\alpha,\alpha}\ast S_{\alpha,\alpha})(t)x+(g_2\ast S_{\alpha,\alpha})(t)x].$
  \item $S_{\alpha,\alpha}(t)x+tS_{\alpha,\alpha}'(t)x=\alpha[A(S_{\alpha,\alpha}\ast S_{\alpha,\alpha})(t)x+S_{\alpha,\alpha}(t)x].$
\end{enumerate}
\end{lemma}
\begin{proof}
The proof of the first assertion follows similarly to the proof of Lemma \ref{Lemma2}. To prove the second one, we integrate (1) to obtain
\begin{equation*}
\int_0^t rS_{\alpha,\alpha}(r)xdr=\alpha(g_1\ast S_{\alpha,1}\ast S_{\alpha,\alpha}(t)x=:\alpha h(t)x.
\end{equation*}
As $\lambda^\alpha(\lambda^\alpha-A)^{-1}=A(\lambda^\alpha-A)^{-1}+I,$ for any $\lambda^\alpha\in\rho(A),$ we have
\begin{equation*}
\hat{h}(\lambda)x=\lambda^{-2}\lambda^\alpha(\lambda^\alpha-A)^{-1}(\lambda^\alpha-A)^{-1}=\lambda^{-2}A(\lambda^\alpha-A)^{-1}(\lambda^\alpha-A)^{-1}+\lambda^{-2}(\lambda^\alpha-A)^{-1},
\end{equation*}
and (2) follows by the uniqueness of the Laplace transform. Finally, by Proposition \ref{Proposition1}, we have
\begin{equation*}
\widehat{S}_{\alpha,1}'(\lambda)x=\lambda \hat{S}_{\alpha,1}(\lambda)x-S_{\alpha,1}(0)x=\lambda^\alpha(\lambda^\alpha-A)^{-1}x-x=A(\lambda^\alpha-A)^{-1}x=A\hat{S}_{\alpha,\alpha}(\lambda)x,
\end{equation*}
and therefore $S_{\alpha,1}'(t)x=AS_{\alpha,\alpha}(t)x.$ To conclude, we notice that (3) is exactly the derivative of (1).
\end{proof}

\section{Determination of $\alpha$ for small times. The sub-diffusion case: $0<\alpha<1.$}\label{Sect3}

Let $0<\alpha<1.$ In this section we determinate the order $\alpha$ of the fractional differential equations $\partial_t^\alpha u(t)=Au(t)$ and $^R\partial_t^\alpha u(t)=Au(t),$ for the Caputo and Riemann-Liouville fractional derivatives. Here $A$ is a given closed linear operator.

Let $x\in X.$ Consider the equation for the Caputo fractional derivative
\begin{equation}\label{eqCaputo-A}
\left\{\begin{array}{lcl}
\partial_t^\alpha u(t) &=& Au(t), \quad t>0\\
u(0)&=&x\\
\end{array}\right.
\end{equation}
where $0<\alpha<1.$ For each $x\in X,$ we denote by $u(t;x)$ the solution to Problem \eqref{eqCaputo-A}. If  $\{S_{\alpha,1}(t)\}_{t\geq 0}$ is the $(\alpha,1)$-resolvent family generated by $A,$ then,
\begin{equation*}\label{Sol-eqCaputo-A}
u(t;x)=S_{\alpha,1}(t)x,
\end{equation*}
see for instance \cite[Chapter 1]{Ba-01}. Next, we define the operators $\varphi_t:X\to X$ and $\psi_t:X\to X,$ respectively, by $\varphi_t(x)=u(t;x)-x$ and $\psi_t(x)=u'(t;x),$ where $u(t;x)$ is the solution of \eqref{eqCaputo-A}. Let $\varphi,\psi:X\to X$ be the operators defined, respectively, by
\begin{equation*}
\varphi(x)=\lim_{t\to 0^+ }t^{-\alpha}\varphi_t(x),\quad \psi(x)=\lim_{t\to 0^+ }t^{1-\alpha}\psi_t(x).
\end{equation*}

Similarly, if $A$ generates the $(\alpha,\alpha)$-resolvent family $\{S_{\alpha,\alpha}(t)\}_{t\geq 0}$ and we consider the Riemann-Liouville in the fractional differential equation
\begin{equation}\label{eqRL-A}
\left\{\begin{array}{lcl}
^R\partial_t^\alpha u(t) &=& Au(t), \quad t>0\\
(g_{1-\alpha}\ast u)(0)&=&x,\\
\end{array}\right.
\end{equation}
then, its solution is given by
\begin{equation*}\label{Sol-eqRL-A}
u(t;x)=S_{\alpha,\alpha}(t)x.
\end{equation*}

Moreover, we define the operators $\tilde{\varphi}_t:X\to X$ and $\tilde{\psi}_t:X\to X,$ respectively, by $\tilde{\varphi}_t(x)=\int_0^t u(s;x)ds=(g_1\ast u)(t;x)$ and $\tilde{\psi}_t(x)=u(t;x),$ where $u(t;x)$ is the solution of \eqref{eqRL-A}. Finally, we define $\tilde{\varphi},\tilde{\psi}:X\to X,$ respectively, by
\begin{equation*}
\tilde{\varphi}(x)=\lim_{t\to 0^+}t^{-\alpha} \tilde{\varphi}_t(x),\quad \tilde{\psi}(x)=\lim_{t\to 0^+} t^{1-\alpha}\tilde{\psi}_t(x).
\end{equation*}

Now, we consider the following problem:
\begin{itemize}
  \item Let $x\in X$ be fixed. Determinate $\alpha\in (0,1)$ in \eqref{eqCaputo-A} and \eqref{eqRL-A} from the observation data $u(t;x)$ for small $t.$
\end{itemize}


The next result gives an answer for the Caputo fractional derivative.
\begin{theorem}\label{Th3.1}
If $A$ generates the $(\alpha,1)$-resolvent family $\{S_{\alpha,1}(t)\}_{t\geq 0}$ and $x\in D(A)\cap D(A^{-1})$ with $A(g_\alpha\ast S_{\alpha,1})(t)x\neq 0$ for $t>0$ small enough, then
\begin{equation*}
\alpha x=\lim_{t\to 0^+} t\psi_t\varphi^{-1}_t(x).
\end{equation*}
\end{theorem}
\begin{proof}
Since $u(t;x)=S_{\alpha,1}(t)x,$ from Proposition \ref{Proposition1} and Lemma \ref{Lemma1}  we can write
\begin{equation}\label{EqCaputo1}
u(t;x)-x=A(g_\alpha\ast S_{\alpha,1})(t)x=AS_{\alpha,\alpha+1}(t)x,
\end{equation}
for any $x\in X,$ and therefore
\begin{equation}\label{EqCaputo0}
t^{-\alpha}\varphi_t(x)=t^{-\alpha}(u(t;x)-x)=\frac{AS_{\alpha,\alpha+1}(t)}{g_{\alpha+1}(t)}\frac{1}{\Gamma(\alpha+1)}x,\quad x\in X.
\end{equation}
By Proposition \ref{Proposition1} we have for any $x\in D(A)$ that
\begin{equation*}\label{EqCaputo2}
\varphi(x)=\lim_{t\to 0^+} t^{-\alpha}\varphi_t(x)=\frac{1}{\Gamma(\alpha+1)}Ax.
\end{equation*}

Now, we claim that $\Gamma(\alpha+1)t^{-\alpha}A^{-1}{\varphi_t}$ is an invertible operator for $t>0$ small enough. In fact, by \eqref{EqCaputo0}, we have
\begin{eqnarray*}
\Gamma(\alpha+1)t^{-\alpha}A^{-1}{\varphi_t}(x)-x=\frac{S_{\alpha,\alpha+1}(t)}{g_{\alpha+1}(t)}x-x.
\end{eqnarray*}
By Proposition \ref{Proposition1} the right hand side in the last identity goes to $0$ as $t\to 0^+.$ Hence, we can take $t>0$ small enough, with $\|\Gamma(\alpha+1)t^{-\alpha}A^{-1}{\varphi_t}(x)-x\|<1$ for all $x\in D(A).$ This implies that $\Gamma(\alpha+1)t^{-\alpha}A^{-1}{\varphi_t}$ is invertible for $t>0$ small enough, and thus,
\begin{equation}\label{EqCaputo4}
\varphi^{-1}(x)=\lim_{t\to 0^+} t^{\alpha}\varphi_t^{-1}(x)=\Gamma(\alpha+1)A^{-1}x,
\end{equation}
for all $x\in D(A)\cap D(A^{-1}).$ On the other hand, for any $\lambda^\alpha\in \rho(A)$ we have by Proposition \ref{Proposition1} that
\begin{equation}\label{EqCaputo5}
\widehat{S}_{\alpha,1}'(\lambda)x=\lambda \hat{S}_{\alpha,1}(\lambda)x-S_{\alpha,1}(0)x=\lambda^\alpha(\lambda^\alpha-A)^{-1}x-x=A(\lambda^\alpha-A)^{-1}x=A\hat{S}_{\alpha,\alpha}(\lambda)x.
\end{equation}
By Proposition \ref{Proposition1} and Lemma \ref{Lemma1} we have for any $x\in D(A)$ that
\begin{equation*}
u'(t;x)=S'_{\alpha,1}(t)x=AS_{\alpha,\alpha}(t)x=g_\alpha(t)Ax+A(g_\alpha\ast S_{\alpha,\alpha})(t)Ax=g_\alpha(t)Ax+AS_{\alpha,2\alpha}(t)Ax.
\end{equation*}
Therefore,
\begin{equation*}
t^{1-\alpha}u'(t;x)=\frac{Ax}{\Gamma(\alpha)}+At^{1-\alpha}S_{\alpha,2\alpha}(t)Ax.
\end{equation*}
The Proposition \ref{Proposition1} implies again that
\begin{equation*}
\lim_{t\to 0^+} t^{1-\alpha}AS_{\alpha,2\alpha}(t)x=\lim_{t\to 0^+}\frac{S_{\alpha,2\alpha}(t)}{g_{2\alpha}(t)}\frac{t^\alpha}{\Gamma(2\alpha)}Ax=0.
\end{equation*}
Hence, we get
\begin{equation}\label{EqCaputo3}
\psi(x)=\lim_{t\to 0^+} t^{1-\alpha}u'(t;x)=\lim_{t\to 0^+} t^{1-\alpha}\psi_t(x)=\frac{Ax}{\Gamma(\alpha)},
\end{equation}
for all $x\in D(A).$ By \eqref{EqCaputo4} and \eqref{EqCaputo3} we obtain
\begin{eqnarray*}
\lim_{t\to 0^+}\|t\psi_t(\varphi_t^{-1}(x))-\psi(\varphi^{-1}(x))\|&=&\lim_{t\to 0^+} \|t^{1-\alpha}\psi_t(t^{\alpha}\varphi_t^{-1}(x)-\varphi^{-1}(x))-\psi(\varphi^{-1}(x))+t^{1-\alpha}\psi_t(\varphi^{-1}(x))\|\\
&\leq& \lim_{t\to 0^+}t^{1-\alpha}\|\psi_t\|\|t^{\alpha}\varphi_t^{-1}(x)-\varphi^{-1}(x)\|+\|(\psi-t^{1-\alpha}\psi_t)(\varphi^{-1}(x))\|\\
&=&0,
\end{eqnarray*}
for all $x\in D(A)\cap D(A^{-1}).$ As
\begin{equation*}
\psi(\varphi^{-1}(x))=\frac{1}{\Gamma(\alpha)}A(\Gamma(\alpha+1)A^{-1}x)=\alpha x,
\end{equation*}
we conclude that
\begin{equation*}
\lim_{t\to 0^+} t\psi_t\varphi^{-1}_t(x)=\alpha x,
\end{equation*}
for all $x\in D(A)\cap D(A^{-1}).$
\end{proof}

Now, we have an answer to the inverse problem for the Riemann-Liouville fractional derivative.

\begin{theorem}
If $A$ generates the $(\alpha,\alpha)$-resolvent family $\{S_{\alpha,\alpha}(t)\}_{t\geq 0}$ and $x\in X$ with $(g_1\ast S_{\alpha,\alpha})(t)x\neq$ for $t>0$ small enough, then
\begin{equation*}
\alpha x=\lim_{t\to 0^+} t\tilde{\psi}_t\tilde{\varphi}_t^{-1}(x).
\end{equation*}
\end{theorem}
\begin{proof}
Let $x\in X.$ As $u(t;x)=S_{\alpha,\alpha}(t)x,$ we have by Proposition \ref{Proposition1} and Lemma \ref{Lemma1} that
\begin{equation*}
t^{1-\alpha}u(t;x)=\frac{x}{\Gamma(\alpha)}+At^{1-\alpha}(g_\alpha\ast S_{\alpha,\alpha})(t)x=\frac{x}{\Gamma(\alpha)}+At^{1-\alpha}S_{\alpha,2\alpha}(t)x.
\end{equation*}
By Proposition \ref{Proposition1} we have
\begin{equation*}
\lim_{t\to 0^+} t^{1-\alpha}AS_{\alpha,2\alpha}(t)x=\lim_{t\to 0^+}\frac{AS_{\alpha,\alpha+1}(t)}{g_{2\alpha}(t)}\frac{t^\alpha}{\Gamma(2\alpha)}x=0.
\end{equation*}
Hence,
\begin{equation}\label{eqRL-0}
\tilde{\psi}(x)=\lim_{t\to 0^+}t^{1-\alpha}\tilde{\psi}_t(x)=\lim_{t\to 0^+} t^{1-\alpha}u(t;x)=\frac{x}{\Gamma(\alpha)}.
\end{equation}
Since $u(t;x)=g_\alpha(t)x+A(g_\alpha \ast S_{\alpha,\alpha})(t)x,$ integrating over $[0,t],$ by Lemma \ref{Lemma1} and the semigroup law for the functions $g_\beta,$ we have that
\begin{equation*}
(g_1\ast u)(t;x)=g_{\alpha+1}(t)x+A(g_1\ast g_\alpha\ast S_{\alpha,\alpha})(t)x=g_{\alpha+1}(t)x+AS_{\alpha,2\alpha+1}(t)x.
\end{equation*}
Therefore,
\begin{equation*}
t^{-\alpha}(g_1\ast u)(t;x)=\frac{x}{\Gamma(\alpha+1)}+At^{-\alpha}S_{\alpha,2\alpha+1}(t)x.
\end{equation*}
By Proposition \ref{Proposition1} we get
\begin{equation*}
\lim_{t\to 0^+} t^{-\alpha}AS_{\alpha,2\alpha+1}(t)x=\lim_{t\to 0^+}\frac{AS_{\alpha,2\alpha+1}(t)}{g_{2\alpha+1}(t)}\frac{t^\alpha}{\Gamma(2\alpha+1)}x=0.
\end{equation*}
Hence,
\begin{equation}\label{EqRL-1}
\tilde{\varphi}(x)=\lim_{t\to 0^+}t^{-\alpha}\tilde{\varphi}_t(x)=\lim_{t\to 0^+}t^{-\alpha}(g_1\ast u)(t;x)=\frac{x}{\Gamma(\alpha+1)}.
\end{equation}
Now, we claim that $t^{-\alpha}\tilde{\varphi_t}$ is an invertible operator for $t>0$ small enough. In fact, as $u(t;x)=S_{\alpha,\alpha}(t)x,$ we have by Proposition \ref{Proposition1} that
\begin{equation*}
\tilde{\varphi}_t(x)=(g_1\ast u)(t;x)=(g_1\ast S_{\alpha,\alpha})(t)x=g_{\alpha+1}(t)x+A(g_\alpha\ast S_{\alpha,\alpha})(t)x,
\end{equation*}
and therefore, the Lemma \ref{Lemma1} implies that
\begin{equation*}
\frac{1}{g_{\alpha+1}(t)}\tilde{\varphi}_t(x)-x=\frac{S_{\alpha,\alpha+1}(t)x}{g_{\alpha+1}(t)}-x=\frac{1}{g_{\alpha+1}(t)}A(g_\alpha\ast S_{\alpha,\alpha})(t)x, \quad t>0.
\end{equation*}
By Proposition \ref{Proposition1} we have $\lim_{t\to 0^+} \frac{S_{\alpha,\alpha+1}(t)x}{g_{\alpha+1}(t)}-x=0,$ and thus we can take $t>0$ small enough, with $\|\frac{1}{g_{\alpha+1}(t)}\tilde{\varphi}_t(x)-x\|<1$ for all $x\in X.$ This implies that $\frac{1}{g_{\alpha+1}(t)}\tilde{\varphi}_t=\Gamma(\alpha+1)t^{-\alpha}\tilde{\varphi}_t$ is an invertible operator, and therefore
$t^{-\alpha}\tilde{\varphi}_t$ is invertible for $t>0$ small enough. By \eqref{EqRL-1} we get
\begin{equation}\label{EqRL-2}
\tilde{\varphi}^{-1}(x)=\lim_{t\to 0^+}t^{\alpha}\tilde{\varphi}^{-1}_t(x)=\Gamma(\alpha+1)x.
\end{equation}

By \eqref{eqRL-0} and \eqref{EqRL-2} have
\begin{eqnarray*}
\lim_{t\to 0^+}\|t\tilde{\psi}_t(\tilde{\varphi}_t^{-1}(x))-\tilde{\psi}(\tilde{\varphi}^{-1}(x))\|&=&\lim_{t\to 0^+} \|t^{1-\alpha}\tilde{\psi}_t(t^{\alpha}\tilde{\varphi}_t^{-1}(x)-\tilde{\varphi}^{-1}(x))-\tilde{\psi}(\tilde{\varphi}^{-1}(x))+t^{1-\alpha}\tilde{\psi}_t(\tilde{\varphi}^{-1}(x))\|\\
&\leq& \lim_{t\to 0^+}t^{1-\alpha}\|\tilde{\psi}_t\|\|t^{\alpha}\tilde{\varphi}_t^{-1}(x)-\tilde{\varphi}^{-1}(x)\|+\|(\tilde{\psi}-t^{1-\alpha}\tilde{\psi}_t)(\tilde{\varphi}^{-1}(x))\|\\
&=&0,
\end{eqnarray*}
for all $x\in X.$ Since $\tilde{\psi}(\tilde{\varphi}^{-1}(x))=\tilde{\psi}(\Gamma(\alpha+1)x)=\frac{\Gamma(\alpha+1)}{\Gamma(\alpha)}x=\alpha x,$ we conclude that
\begin{equation*}
\lim_{t\to 0^+} t\tilde{\psi}_t\tilde{\varphi}_t^{-1}(x)=\alpha x.
\end{equation*}
\end{proof}

\section{Determination of $\alpha$ for small times. The super-diffusion case: $1<\alpha<2.$}\label{Sect4}
In this section we determinate the order $\alpha$ (where $1<\alpha<2$) of the fractional differential equations $\partial_t^\alpha u(t)=Au(t)$ and $^R\partial_t^\alpha u(t)=Au(t),$ for the Caputo and Riemann-Liouville fractional derivatives.

We first consider the  equation for the Caputo fractional derivative
\begin{equation}\label{eqCaputo-B}
\left\{\begin{array}{lcl}
\partial_t^\alpha u(t) &=& Au(t), \quad t\geq 0\\
u(0)&=&x\\
u'(0)&=&y.
\end{array}\right.
\end{equation}
For each $x,y\in X,$ we denote by $u(t;x,y)$ the solution to problem \eqref{eqCaputo-B}. If  $\{S_{\alpha,1}(t)\}_{t\geq 0}$ is the $(\alpha,1)$-resolvent family generated by $A,$ then,
\begin{equation*}\label{Sol-eqCaputo-B}
u(t;x,y)=S_{\alpha,1}(t)x+(g_1\ast S_{\alpha,1})(t)y,
\end{equation*}
see for instance \cite{Pon-20}. By Lemma \ref{Lemma1} we can write
\begin{equation*}\label{Sol-eqCaputo-B}
u(t;x,y)=S_{\alpha,1}(t)x+S_{\alpha,2}(t)y.
\end{equation*}

Now, we define the operators $\varphi_t:X\to X$ and $\psi_t:X\to X,$ respectively, by $\varphi_t(x)=u(t;x,y)-x-ty$ and $\psi_t(x)=u''(t;x,y),$ where $u(t;x,y)$ is the solution of \eqref{eqCaputo-B} for any fixed $y\in X.$ Moreover, we define $\varphi,\psi:X\to X,$ respectively, by
\begin{equation*}
\varphi(x)=\lim_{t\to 0^+} t^{-\alpha}\varphi_t(x),\quad \psi(x)=\lim_{t\to 0^+} t^{2-\alpha}\psi_t(x),\quad x\in X.
\end{equation*}

Now, we consider the Riemann-Liouville fractional derivative. If $u(t;x,y)$ denotes the solution to the problem
\begin{equation}\label{eqRL-B}
\left\{\begin{array}{lcl}
\qquad \,^R\partial_t^\alpha u(t) &=& Au(t), \quad t\geq 0\\
(g_{2-\alpha}\ast u)(0)&=&x\\
(g_{2-\alpha}\ast u)'(0)&=&y,
\end{array}\right.
\end{equation}
and $A$ generates the $(\alpha,\alpha-1)$-resolvent family $\{S_{\alpha,\alpha-1}(t)\}_{t\geq 0},$ then (by \cite{Pon-20}), we have
\begin{equation}\label{Sol1-eqRL-B}
u(t;x,y)=S_{\alpha,\alpha-1}(t)x+(g_1\ast S_{\alpha,\alpha-1})(t)y.
\end{equation}
Moreover, by Lemma \ref{Lemma1} we can write
\begin{equation}\label{Sol-eqRL-B}
u(t;x,y)=S_{\alpha,\alpha-1}(t)x+S_{\alpha,\alpha}(t)y.
\end{equation}

Moreover, for any fixed $y\in X,$ we define the operators $\tilde{\varphi}_t:X\to X$ and $\tilde{\psi}_t:X\to X,$ respectively, by $\tilde{\varphi}_t(x)=\int_0^t u(s;x,y)ds=(g_1\ast u)(t;x,y)$ and $\tilde{\psi}_t(x)=\int_0^t (t-s)u(s;x,y)ds=(g_2\ast u)(t;x,y),$ where $u(t;x,y)$ is the solution of \eqref{eqRL-B}. Finally, we define $\tilde{\varphi},\tilde{\psi}:X\to X,$ respectively, by
\begin{equation*}
\tilde{\varphi}(x)=\lim_{t\to 0^+} t^{1-\alpha}\tilde{\varphi}_t(x),\quad \tilde{\psi}(x)=\lim_{t\to 0^+} t^{-\alpha}\tilde{\psi}_t(x),\quad x\in X.
\end{equation*}

Next, we consider the following problem:
\begin{itemize}
  \item Let $x\in X$ be fixed. Determinate $\alpha\in (1,2)$ in \eqref{eqCaputo-B} and \eqref{eqRL-B} from the observation data $u(t;x,y)$ for small $t.$
\end{itemize}

For the Caputo fractional derivative \eqref{eqCaputo-B} we have the following result.

\begin{theorem}\label{Th4.3}
Let $y\in X.$ If $A$ generates the $(\alpha,1)$-resolvent family $\{S_{\alpha,1}(t)\}_{t\geq 0}$ and $x\in D(A)\cap D(A^{-1})$ with $S_{\alpha,1}(t)x+(g_1\ast S_{\alpha,1})(t)y\neq 0$ for $t>0$ small enough, then
\begin{equation*}
\alpha(\alpha-1)x=\lim_{t\to 0^+} t^2\psi_t\varphi_t^{-1}(x).
\end{equation*}
\end{theorem}
\begin{proof}
Let $x\in D(A)\cap D(A^{-1}).$ As $u(t;x,y)=S_{\alpha,1}(t)x+S_{\alpha,2}(t)y,$ by Proposition \ref{Proposition1} and Lemma \ref{Lemma1} we can write
\begin{equation}\label{EqCaputo6}
u(t;x,y)=x+A(g_\alpha \ast S_{\alpha,1})(t)x+ty+A(g_\alpha\ast S_{\alpha,2})(t)y=x+AS_{\alpha,\alpha+1}(t)x+ty+AS_{\alpha,\alpha+2}(t)y.
\end{equation}
The Proposition \ref{Proposition1} implies
\begin{equation*}
\lim_{t\to 0^+} t^{-\alpha}AS_{\alpha,\alpha+1}(t)x=\lim_{t\to 0^+}\frac{AS_{\alpha,\alpha+1}(t)}{g_{\alpha+1}(t)}\frac{1}{\Gamma(\alpha+1)}Ax=\frac{1}{\Gamma(\alpha+1)}Ax
\end{equation*}
and
\begin{equation*}
\lim_{t\to 0^+} t^{-\alpha}AS_{\alpha,\alpha+2}(t)x=\lim_{t\to 0^+}\frac{AS_{\alpha,\alpha+2}(t)}{g_{\alpha+2}(t)}\frac{t}{\Gamma(\alpha+2)}Ax=0.
\end{equation*}
Hence,
\begin{equation*}
\lim_{t\to 0^+} t^{-\alpha}(u(t;x,y)-x-ty)=\frac{1}{\Gamma(\alpha+1)}Ax,
\end{equation*}
that is,
\begin{equation}\label{EqCaputo7}
\varphi(x)=\lim_{t\to 0^+} t^{-\alpha}\varphi_t(x)=\frac{1}{\Gamma(\alpha+1)}Ax,\quad x\in D(A), y\in X.
\end{equation}
By \eqref{EqCaputo5}, $S_{\alpha,1}'(t)=AS_{\alpha,\alpha}(t)$ and $S_{\alpha,2}'(t)=S_{\alpha,1}(t),$ therefore
\begin{equation*}
u'(t;x,y)=S_{\alpha,1}'(t)x+S_{\alpha,2}'(t)y=AS_{\alpha,\alpha}(t)x+S_{\alpha,1}(t)y.
\end{equation*}
As $\alpha>1,$ we have for any $\lambda^\alpha\in \rho(A),$
\begin{equation}\label{eqCaputo5b}
\widehat{S}_{\alpha,\alpha}'(\lambda)x=\lambda \hat{S}_{\alpha,\alpha}(\lambda)x-S_{\alpha,\alpha}(0)x=\lambda^{\alpha-(\alpha-1)}(\lambda^\alpha-A)^{-1}x=\hat{S}_{\alpha,\alpha-1}(\lambda)x.
\end{equation}
Hence,
\begin{equation*}
t^{2-\alpha}u''(t;x,y)=t^{2-\alpha}AS_{\alpha,\alpha}'(t)x+t^{2-\alpha}S_{\alpha,1}'(t)y=t^{2-\alpha}AS_{\alpha,\alpha-1}(t)x+t^{2-\alpha}AS_{\alpha,\alpha}(t)y.
\end{equation*}
By Proposition \ref{Proposition1} we have
\begin{equation*}
\lim_{t\to 0^+} t^{2-\alpha}AS_{\alpha,\alpha-1}(t)x=\lim_{t\to 0^+}\frac{AS_{\alpha,\alpha-1}(t)}{g_{\alpha-1}(t)}\frac{1}{\Gamma(\alpha-1)}x=\frac{Ax}{\Gamma(\alpha-1)}
\end{equation*}
and
\begin{equation*}
\lim_{t\to 0^+} t^{2-\alpha}AS_{\alpha,\alpha}(t)x=\lim_{t\to 0^+}\frac{AS_{\alpha,\alpha}(t)}{g_{\alpha}(t)}\frac{t}{\Gamma(\alpha)}x=0.
\end{equation*}
Therefore,
\begin{equation*}
\lim_{t\to 0^+} t^{2-\alpha}u''(t;x,y)=\frac{Ax}{\Gamma(\alpha-1)},
\end{equation*}
that is,
\begin{equation}\label{EqCaputo8}
\psi(x)=\lim_{t\to 0^+} t^{2-\alpha}\psi_t(x)=\frac{Ax}{\Gamma(\alpha-1)}.
\end{equation}

Now, we claim that $\Gamma(\alpha+1)t^{-\alpha}\varphi_t$ is invertible for $t>0$ small enough. In fact, by \eqref{EqCaputo6} we can write
\begin{eqnarray*}
\Gamma(\alpha+1)t^{-\alpha}A^{-1}{\varphi_t}(x)-x=\frac{S_{\alpha,\alpha+1}(t)}{g_{\alpha+1}(t)}x-x+\frac{S_{\alpha,\alpha+2}(t)}{g_{\alpha+2}(t)}\frac{g_{\alpha+2}(t)}{g_{\alpha+1}(t)}y.
\end{eqnarray*}
The Proposition \ref{Proposition1} implies that the right hand side in the last identity goes to $0$ as $t\to 0^+.$ Hence, we take $t>0$ small enough, with $\|\Gamma(\alpha+1)t^{-\alpha}A^{-1}{\varphi_t}(x)-x\|<1$ for all $x\in D(A).$ This implies that $\Gamma(\alpha+1)t^{-\alpha}A^{-1}{\varphi_t}$ is invertible for $t>0$ small enough, and thus,
\begin{equation}\label{EqCaputo9}
\varphi^{-1}(x)=\lim_{t\to 0^+}(t^{-\alpha} \varphi_t)^{-1}(x)=\Gamma(\alpha+1)A^{-1}x,
\end{equation}
for all $x\in D(A)\cap D(A^{-1}).$ Now, by \eqref{EqCaputo8} and \eqref{EqCaputo9} we get
\begin{eqnarray*}
\lim_{t\to 0^+}\|t^2\psi_t(\varphi_t^{-1}(x))&-&\psi(\varphi^{-1}(x))\|\\
&=&\lim_{t\to 0^+} \|t^{2-\alpha}\psi_t(t^{\alpha}\varphi_t^{-1}(x)-\varphi^{-1}(x))-\psi(\varphi^{-1}(x))+t^{2-\alpha}\psi_t(\varphi^{-1}(x))\|\\
&\leq& \lim_{t\to 0^+}t^{2-\alpha}\|\psi_t\|\|t^{\alpha}\varphi_t^{-1}(x)-\varphi^{-1}(x)\|+\|(\psi-t^{2-\alpha}\psi_t)(\varphi^{-1}(x))\|\\
&=&0,
\end{eqnarray*}
for all $x\in D(A)\cap D(A^{-1}).$ As $\psi(\varphi^{-1}(x))=\psi(\Gamma(\alpha+1)A^{-1}x)=\frac{1}{\Gamma(\alpha-1)}\Gamma(\alpha+1)x=\alpha(\alpha-1)x,$ we conclude that
\begin{equation*}
\lim_{t\to 0^+} t^2\psi_t \varphi_t^{-1}(x)=\alpha(\alpha-1)x.
\end{equation*}
\end{proof}

Now, we consider the problem \eqref{eqRL-B} for the Riemann-Liouville fractional derivative.

\begin{theorem}
Let $y\in X.$ If $A$ generates the $(\alpha,\alpha-1)$-resolvent family $\{S_{\alpha,\alpha-1}(t)\}_{t\geq 0}$ and for $x\in X$ we have $S_{\alpha,\alpha-1}(t)x+(g_1\ast S_{\alpha,\alpha-1})(t)y\neq 0$ for $t>0$ small enough, then
\begin{equation*}
\alpha x=\lim_{t\to 0^+} t^2\tilde{\varphi}_t\tilde{\psi}_t^{-1}(x).
\end{equation*}
\end{theorem}
\begin{proof}
Let $x\in X.$ As $u(t;x,y)=S_{\alpha,\alpha-1}(t)x+S_{\alpha,\alpha}(t)y,$ by Proposition \ref{Proposition1}, Lemma \ref{Lemma1}, and the semigroup law for the functions $g_\beta,$ we obtain
\begin{eqnarray*}
(g_1\ast u)(t;x,y)&=&\int_0^t u(s;x,y)ds\\
&=&g_\alpha(t)x+A(g_{\alpha+1}\ast S_{\alpha,\alpha-1})(t)x+g_{\alpha+1}(t)y+A(g_{\alpha+1}\ast S_{\alpha,\alpha})(t)y\\
&=&g_\alpha(t)x+AS_{\alpha,2\alpha}(t)x+g_{\alpha+1}(t)y+AS_{\alpha,2\alpha+1}(t)y.
\end{eqnarray*}
Hence,
\begin{equation*}
t^{1-\alpha}(g_1\ast u)(t;x,y)=\frac{1}{\Gamma(\alpha)}x+t^{1-\alpha}AS_{\alpha,2\alpha}(t)x+\frac{ty}{\Gamma(\alpha+1)}+t^{1-\alpha}AS_{\alpha,2\alpha+1}(t)y.
\end{equation*}
The Proposition \ref{Proposition1} implies that
\begin{equation*}
\lim_{t\to 0^+} t^{1-\alpha}AS_{\alpha,2\alpha}(t)x=\lim_{t\to 0^+}\frac{AS_{\alpha,2\alpha}(t)}{g_{2\alpha}(t)}\frac{t^{\alpha}}{\Gamma(2\alpha)}x=0
\end{equation*}
and
\begin{equation*}
\lim_{t\to 0^+} t^{1-\alpha}AS_{\alpha,2\alpha+1}(t)y=\lim_{t\to 0^+}\frac{AS_{\alpha,2\alpha+1}(t)}{g_{2\alpha+1}(t)}\frac{t^{\alpha+1}}{\Gamma(2\alpha+1)}y=0.
\end{equation*}
We obtain
\begin{equation*}
\lim_{t\to 0^+}t^{1-\alpha}(g_1\ast u)(t;x,y)=\frac{1}{\Gamma(\alpha)}x,
\end{equation*}
which means that
\begin{equation}\label{EqRL-5}
\tilde{\varphi}(x)=\lim_{t\to 0^+}t^{1-\alpha}\tilde{\varphi}_t(x)=\frac{1}{\Gamma(\alpha)}x,
\end{equation}
for all $x,y\in X.$

Now, we integrate \eqref{Sol-eqRL-B} twice to obtain, by Proposition \ref{Proposition1} and Lemma \ref{Lemma1}, that
\begin{eqnarray}\label{eqRL-3}
(g_2\ast u)(t;x,y)&=&\int_0^t (t-s)u(s;x,y)ds \nonumber \\
&=&g_{\alpha+1}(t)x+A(g_{\alpha+2}\ast S_{\alpha,\alpha-1})(t)x+g_{\alpha+2}(t)y+A(g_{\alpha+2}\ast S_{\alpha,\alpha})(t)y\nonumber\\
&=&g_{\alpha+1}(t)x+AS_{\alpha,2\alpha+1}(t)x+g_{\alpha+2}(t)y+AS_{\alpha,2\alpha+2}(t)y.
\end{eqnarray}
Thus
\begin{equation*}
t^{-\alpha}(g_2\ast u)(t;x,y)=\frac{1}{\Gamma(\alpha+1)}x+t^{-\alpha}AS_{\alpha,2\alpha+1}(t)x+\frac{ty}{\Gamma(\alpha+2)}+t^{-\alpha}AS_{\alpha,2\alpha+2}(t)y.
\end{equation*}
As
\begin{equation*}
\lim_{t\to 0^+} t^{-\alpha}AS_{\alpha,2\alpha+1}(t)x=\lim_{t\to 0^+}\frac{AS_{\alpha,2\alpha+1}(t)}{g_{2\alpha+1}(t)}\frac{t^{\alpha}}{\Gamma(2\alpha+1)}x=0
\end{equation*}
and
\begin{equation*}
\lim_{t\to 0^+} t^{-\alpha}AS_{\alpha,2\alpha+2}(t)y=\lim_{t\to 0^+}\frac{AS_{\alpha,2\alpha+2}(t)}{g_{2\alpha+2}(t)}\frac{t^{\alpha+1}}{\Gamma(2\alpha+2)}y=0,
\end{equation*}
(see Proposition \ref{Proposition1}) we conclude that
\begin{equation}\label{EqRL-6}
\tilde{\psi}(x)=\lim_{t\to 0^+}t^{-\alpha}\tilde{\psi}_t(x)=\lim_{t\to 0^+}t^{-\alpha}(g_2\ast u)(t;x,y)=\frac{1}{\Gamma(\alpha+1)}x.
\end{equation}
Now, we will see that $t^{-\alpha}\tilde{\psi}_t$ is invertible for $t>0$ small enough and any $y\in X$ being fixed. In fact, by \eqref{eqRL-3}, we obtain
\begin{eqnarray*}
\frac{1}{g_{\alpha+1}(t)}\tilde{\psi}_t(x)-x&=&\frac{S_{\alpha,2\alpha+1}(t)}{g_{\alpha+1}(t)}x+\frac{g_{\alpha+2}(t)}{g_{\alpha+1}(t)}y+\frac{AS_{\alpha,2\alpha+2}(t)}{g_{\alpha+1}(t)}y\\
&=&\frac{S_{\alpha,2\alpha+1}(t)}{g_{2\alpha+1}(t)}\frac{g_{2\alpha+1}(t)}{g_{\alpha+1}(t)}x+\frac{g_{\alpha+2}(t)}{g_{\alpha+1}(t)}y+\frac{AS_{\alpha,2\alpha+2}(t)}{g_{2\alpha+2}(t)}\frac{g_{2\alpha+2}(t)}{g_{\alpha+1}(t)}y \quad t>0.
\end{eqnarray*}
By Proposition \ref{Proposition1}, the right hand side in the last equality goes to $0$ as $t\to 0^+,$ and therefore, we can choose $t>0$ small enough such that $\|\frac{1}{g_{\alpha+1}(t)}\tilde{\psi}_t(x)-x\|<1$ for all $x\in X.$ This implies that $\frac{1}{g_{\alpha+1}(t)}\tilde{\psi}_t$ is an invertible operator, and therefore
$t^{-\alpha}\tilde{\psi}_t$ is invertible for $t>0$ small enough. By \eqref{EqRL-6} we obtain
\begin{equation*}
\tilde{\psi}^{-1}(x)\lim_{t\to 0^+}(t^{-\alpha}\tilde{\psi}_t)^{-1}(x)=\Gamma(\alpha+1)x.
\end{equation*}
As
\begin{equation*}
\tilde{\varphi}(\tilde{\psi}^{-1}(x))=\tilde{\varphi}(\Gamma(\alpha+1)x)=\frac{1}{\Gamma(\alpha)}\Gamma(\alpha+1)x=\alpha x, \end{equation*}
the conclusion follows as in the proof of Theorem \ref{Th4.3}.
\end{proof}

\section{Determination of $\alpha$ for a fixed time $T.$ The sub-diffusion case: $0<\alpha<1.$ }\label{Sect5}

In this section we consider the problem of finding the order $\alpha\in(0,1)$ for a fixed time $T>0$ in the fractional  problems \eqref{eqCaputo-A} and \eqref{eqRL-A}. We first consider the problem for the Caputo fractional derivative. Assume that $A$ is the generator of the resolvent family $\{S_{\alpha,1}(t)\}_{t\geq 0}.$ Let $\varphi_t:X\to X$ be the operator defined by $\varphi_t(x):=(S_{\alpha,1}\ast u)(t;x)-(g_1\ast u)(t;x),$ where $u(t;x)$ is the solution to Problem \eqref{eqCaputo-A}.

\begin{theorem}\label{Th5.1}
If $A$ generates the $(\alpha,1)$-resolvent family $\{S_{\alpha,1}(t)\}_{t\geq 0},$ $x\in X$ and $T>0$ are fixed, then the order $\alpha$ verifies
\begin{equation*}
Tu(T;x)-(g_1\ast u)(T;x)=\alpha \varphi_T(x).
\end{equation*}
\end{theorem}
\begin{proof}
Let $x\in X$ and $T>0.$ By (2) in Lemma \ref{Lemma2} we have
\begin{equation}\label{EqCaputo10}
tS_{\alpha,1}(t)x-(g_1\ast S_{\alpha,1})(t)x=\alpha[(S_{\alpha,1}\ast S_{\alpha,1})(t)x-(g_1\ast S_{\alpha,1})(t)x]
\end{equation}
for all $t\geq 0.$  As $u(t;x)=S_{\alpha,1}(t)x$ is the solution to Problem \eqref{eqCaputo-A}, we have
\begin{eqnarray*}
\varphi_t(x)&=&\int_{0}^t S_{\alpha,1}(t-r)u(r;x)dr-\int_0^t u(r;x)dr\\
&=&(S_{\alpha,1}\ast S_{\alpha,1})(t)x-(g_1\ast S_{\alpha,1})(t)x.
\end{eqnarray*}
Therefore, \eqref{EqCaputo10} can be written as
\begin{equation*}
tS_{\alpha,1}(t)x-(g_1\ast S_{\alpha,1})(t)x=\alpha \varphi_t(x),
\end{equation*}
for any $t>0$ and $x\in X.$ We conclude that
\begin{equation*}
Tu(T;x)-(g_1\ast u)(T;x)=\alpha \varphi_T(x).
\end{equation*}
\end{proof}

\begin{remark}\label{Rem1}
We notice that if $u(t;x)$ is real valued, then to find $\alpha,$ we only to {\em divide} by $\varphi_T$ in Theorem \ref{Th5.1} to obtain
\begin{equation*}
\alpha=\frac{Tu(T;x)-(g_1\ast u)(T;x)}{\varphi_T(x)}=\frac{Tu(T;x)-(g_1\ast u)(T;x)}{(S_{\alpha,1}\ast u)(T;x)-(g_1\ast u)(T;x)},
\end{equation*}
that is, we need to know the data: The solution $u(T;x),$ its integral $(g_1\ast u)(T;x)$ and the convolution $(S_{\alpha,1}\ast u)(T;x)$ for a fixed $x\in X$ and a time $T>0.$
\end{remark}


On the other hand, we notice that by Lemma \ref{Lemma2} we have
\begin{equation}\label{EqCaputo11}
tS_{\alpha,1}'(t)x=\alpha(S_{\alpha,1}'\ast S_{\alpha,1})(t)x,
\end{equation}
for any $t>0$ and $x\in X.$ As $u(t;x)=S_{\alpha,1}(t)x$ is the solution to Problem \eqref{eqCaputo-A}, if $F_t:X\to X$ is the operator defined by $F_t(x):=(S_{\alpha,1}'\ast S_{\alpha,1})(t)x,$ then the order $\alpha$ also verifies
\begin{equation*}
F_t(x)=\int_0^t S_{\alpha,1}'(r)u(r;x)dr=(S_{\alpha,1}'\ast u)(t;x),
\end{equation*}
and, as in Theorem \ref{Th5.1}, we obtain
\begin{equation*}
tu'(t;x)=\alpha F_t(x)
\end{equation*}
for any $t>0$ and $x\in X.$ Therefore, by \eqref{EqCaputo11} we have the following result.
\begin{theorem}\label{Th5.2}
If $A$ generates the $(\alpha,1)$-resolvent family $\{S_{\alpha,1}(t)\}_{t\geq 0},$ $x\in X$ and $T>0,$ then
\begin{equation*}
Tu'(T;x)=\alpha F_T(x).
\end{equation*}
\end{theorem}

Now, we consider the Problem \eqref{eqRL-A} for the Riemann-Liouville fractional derivative. Assume that $A$ is the generator of $\{S_{\alpha,\alpha}(t)\}_{t\geq 0}.$ Let $\psi_t:X\to X$ be the operator defined by $\psi_t(x):=A(g_2\ast S_{\alpha,\alpha}\ast S_{\alpha,\alpha})(t)x+(g_2\ast S_{\alpha,\alpha})(t)x=A(g_2\ast S_{\alpha,\alpha}\ast u)(t;x)ds+(g_2\ast u)(t;x),$ where $u(t;x)$ is the solution to Problem \eqref{eqRL-A}.

\begin{theorem}\label{Th5.3}
If $A$ generates the $(\alpha,\alpha)$-resolvent family $\{S_{\alpha,\alpha}(t)\}_{t\geq 0}$ and $x\in X,$ $T>0,$ then the order $\alpha$ verifies
\begin{equation*}
\int_0^T ru(r;x)dr=\alpha \psi_T(x).
\end{equation*}
\end{theorem}
\begin{proof}
Let $t>0$ and $x\in X.$ As $u(t;x)=S_{\alpha,\alpha}(t)x$ is the solution to \eqref{eqRL-A}, by Lemma \ref{Lemma3} we have
\begin{equation*}\label{eqRL-4}
\int_0^t rS_{\alpha,\alpha}(r)xdr=\alpha[A(g_2\ast S_{\alpha,\alpha}\ast S_{\alpha,\alpha})(t)x+(g_2\ast S_{\alpha,\alpha})(t)x]=\alpha\psi_t(x),
\end{equation*}
for any $t>0$ and $x\in X.$ We conclude that
\begin{equation*}
\int_0^T ru(r;x)dr=\alpha \psi_T(x),
\end{equation*}
for any $x\in X.$
\end{proof}

Finally, by Lemma \ref{Lemma3}, we notice that $tS_{\alpha,\alpha}(t)x=\alpha(S_{\alpha,1}\ast S_{\alpha,\alpha})(t)x,$ for all $t\geq 0$ and $x\in X.$ As $\lambda^{\alpha}(\lambda^{\alpha}-A)^{-1}=A(\lambda^{\alpha}-A)^{-1}+I,$ we get
\begin{equation*}
\mathcal{L}((S_{\alpha,1}\ast S_{\alpha,\alpha}))(\lambda)x=\lambda^{\alpha-1}(\lambda^{\alpha}-A)^{-1}(\lambda^{\alpha}-A)^{-1}x=\frac{1}{\lambda}A(\lambda^{\alpha}-A)^{-1}(\lambda^{\alpha}-A)^{-1}x+\frac{1}{\lambda}(\lambda^{\alpha}-A)^{-1}x,
\end{equation*}
which implies that
\begin{equation*}
tS_{\alpha,\alpha}(t)x=\alpha [(g_1\ast AS_{\alpha,\alpha}\ast S_{\alpha,\alpha})(t)x+(g_1\ast S_{\alpha,\alpha})(t)x],
\end{equation*}
for all $t\geq 0,$ $x\in X.$ Since $u(t;x)=S_{\alpha,\alpha}(t)x$ is the solution to Problem \eqref{eqRL-A}, we have that if $G_t(x):=A(g_1\ast S_{\alpha,\alpha}\ast S_{\alpha,\alpha})(t)x+(g_1\ast S_{\alpha,\alpha})(t)x,$ then
\begin{equation*}
tu(t;x)=\alpha G_t(x)=\alpha[A(g_1\ast S_{\alpha,\alpha}\ast u)(t;x)+(g_1\ast u)(t;x)]
\end{equation*}
Therefore, we have the following result.

\begin{theorem}\label{Th5.4}
If $A$ generates the $(\alpha,\alpha)$-resolvent family $\{S_{\alpha,\alpha}(t)\}_{t\geq 0}$ and $x\in X,$ $T>0,$ then the order $\alpha$ verifies
\begin{equation*}
Tu(T;x)=\alpha G_T(x).
\end{equation*}
\end{theorem}

\section{Determination of $\alpha$ for a fixed time $T.$ The super-diffusion case: $1<\alpha<2.$ }\label{Sect6}

In this section we find the order $\alpha\in(1,2)$ for a fixed time $T>0$ in the fractional problems \eqref{eqCaputo-B} and \eqref{eqRL-B}. We first consider the problem \eqref{eqCaputo-B}. Assume that $A$ is the generator of the resolvent family $\{S_{\alpha,1}(t)\}_{t\geq 0}.$ For a given $y\in X,$ let $\varphi_t:X\to X$ be the operator defined by $\varphi_t(x):=(S_{\alpha,1}\ast u)(t;x,y)-(g_1\ast u)(t;x,y),$ where $u(t;x,y)$ is the solution to Problem \eqref{eqCaputo-B}.

\begin{theorem}\label{Th6.1}
If $A$ generates the $(\alpha,1)$-resolvent family $\{S_{\alpha,1}(t)\}_{t\geq 0},$ $x,y\in X$ and $T>0,$ then the order $\alpha$ verifies
\begin{equation*}
Tu(T;x,y)-(g_1\ast u)(T;x,y)-(g_2\ast S_{\alpha,1})(T)y=\alpha \varphi_T(x).
\end{equation*}
\end{theorem}
\begin{proof}
Let $x\in X$ and $T>0.$ We first notice that if $h(t)=t(g_1\ast S_{\alpha,1})(t),$ then for any $\lambda^\alpha\in\rho(A)$ we have
\begin{equation*}
\hat{h}(\lambda)x=-\frac{d}{d\lambda}(\mathcal{L}(g_1\ast S_{\alpha,1}))(\lambda)x=-\frac{d}{d\lambda}\left(\lambda^{\alpha-2}(\lambda^\alpha-A)^{-1}x\right)=-(\alpha-2)\lambda^{-2}\hat{S}_{\alpha,1}(\lambda)x+\alpha\lambda^{-1}\hat{S}_{\alpha,1}(\lambda)\hat{S}_{\alpha,1}(\lambda)x.
\end{equation*}
This means that
\begin{equation}\label{EqCaputo12}
t(g_1\ast S_{\alpha,1})(t)x=-(\alpha-2)(g_2\ast S_{\alpha,1})(t)x+\alpha (g_1\ast S_{\alpha,1}\ast S_{\alpha,1})(t)x.
\end{equation}
for all $t\geq 0$ and $x\in X.$ Moreover, by Lemma \ref{Lemma1}, \ref{Lemma2} and \eqref{EqCaputo12} we have
\begin{eqnarray*}
tu(t;x,y)&=&tS_{\alpha,1}(t)x+t(g_1\ast S_{\alpha,1})(t)y\\
&=&\alpha[(S_{\alpha,1}\ast S_{\alpha,1})(t)x+(g_1\ast S_{\alpha,1}\ast S_{\alpha,1})(t)y-(g_1\ast S_{\alpha,1})(t)x-(g_2\ast S_{\alpha,1})(t)y]\\
&&+(g_1\ast S_{\alpha,1})(t)x+2(g_2\ast S_{\alpha,1})(t)y.
\end{eqnarray*}
Now, by Lemma \ref{Lemma1}, we get
\begin{equation*}
tu(t;x,y)-[(g_1\ast S_{\alpha,1})(t)x+(g_1\ast S_{\alpha,2})(t)y]-(g_2\ast S_{\alpha,1})(t)y=\alpha[(S_{\alpha,1}\ast u)(t;x,y)-(g_1\ast u)(t;x,y)],
\end{equation*}
that is,
\begin{equation*}
tu(t;x,y)-(g_1\ast u)(t;x,y)-(g_2\ast S_{\alpha,1})(t)y=\alpha[(S_{\alpha,1}\ast u)(t;x,y)-(g_1\ast u)(t;x,y)].
\end{equation*}
\end{proof}

Finally, we consider Problem \eqref{eqRL-B}. Assume that $A$ is the generator of the resolvent family $\{S_{\alpha,\alpha-1}(t)\}_{t\geq 0}.$ By \eqref{Sol1-eqRL-B}, the solution to \eqref{eqRL-B} is given by
$u(t;x,y)=S_{\alpha,\alpha-1}(t)x+(g_1\ast S_{\alpha,\alpha-1})(t)y.$ For a fixed $y\in X,$ we define $\tilde{\varphi}_t:X\to X$ be the operator defined by $\tilde{\varphi}_t(x):=(g_1\ast u)(t;x,y)+A(g_2\ast S_{\alpha,\alpha-1}\ast u)(t;x,y),$ where $u(t;x,y)$ is the solution to Problem \eqref{eqRL-B}.

\begin{theorem}\label{Th6.2}
If $A$ generates the $(\alpha,\alpha-1)$-resolvent family $\{S_{\alpha,\alpha-1}(t)\}_{t\geq 0},$ $x,y\in X$ and $T>0,$ then the order $\alpha$ verifies
\begin{equation*}
Tu(T;x,y)+(g_1\ast S_{\alpha,\alpha-1})(T)x=\alpha \tilde{\varphi}_T(x).
\end{equation*}
\end{theorem}
\begin{proof}
Let $t>0$ and $x,y\in X.$ By \eqref{Sol-eqRL-B} we have
\begin{equation*}
tu(t;x,y)=tS_{\alpha,\alpha-1}(t)x+tS_{\alpha,\alpha}(t)y.
\end{equation*}
By Lemma \ref{Lemma1} we have $S_{\alpha,\alpha}(t)=(g_1\ast S_{\alpha,\alpha-1})(t)$ and thus $S_{\alpha,\alpha}'(t)=S_{\alpha,\alpha-1}(t).$ Hence
\begin{equation*}
tu(t;x,y)=tS_{\alpha,\alpha}'(t)x+tS_{\alpha,\alpha}(t)y.
\end{equation*}
As in the proof of Lemma \ref{Lemma3}, it is easy to see that
\begin{equation}\label{EqRL-7}
tS_{\alpha,\alpha}(t)=\alpha (S_{\alpha,1}\ast S_{\alpha,\alpha})(t),
\end{equation}
for any $t\geq 0.$ As $S_{\alpha,\alpha}'(t)=S_{\alpha,\alpha-1}(t)$ and for $\alpha>1,$ $S_{\alpha,\alpha}(0)=0,$ we get
\begin{equation}\label{EqRL-8}
S_{\alpha,\alpha}(t)+tS_{\alpha,\alpha}'(t)=\alpha (S_{\alpha,1}\ast S_{\alpha,\alpha}')(t)+\alpha S_{\alpha,1}(t)S_{\alpha,\alpha}'(0)=\alpha (S_{\alpha,1}\ast S_{\alpha,\alpha-1})(t),
\end{equation}
for any $t\geq 0.$  By \eqref{EqRL-7} and \eqref{EqRL-8} we have
\begin{eqnarray*}
tu(t;x,y)&=&tS_{\alpha,\alpha}'(t)x+tS_{\alpha,\alpha}(t)y\\
&=&\alpha[(S_{\alpha,1}\ast S_{\alpha,\alpha-1})(t)x+(S_{\alpha,1}\ast S_{\alpha,\alpha})(t)y]-S_{\alpha,\alpha}(t)x\\
&=&\alpha\int_0^t S_{\alpha,1}(t-s)[S_{\alpha,\alpha-1}(s)x+S_{\alpha,\alpha}(s)y]ds-S_{\alpha,\alpha}(t)x\\
&=&\alpha\int_0^t S_{\alpha,1}(t-s)u(s;x,y)ds-S_{\alpha,\alpha}(t)x.
\end{eqnarray*}
As $S_{\alpha,1}'(t)=AS_{\alpha,\alpha}(t)$ and $S_{\alpha,1}(0)=I,$ integrating by parts, we obtain
\begin{eqnarray*}
\int_0^t S_{\alpha,1}(t-s)u(s;x,y)ds&=&S_{\alpha,1}(t-s)(g_1\ast u)(s;x,y)dr \Big|_{s=0}^{s=t}+\int_0^t AS_{\alpha,\alpha}(t-s)(g_1\ast u)(s;x,y)ds\\
&=&(g_1\ast u)(t;x,y)+A(g_1\ast S_{\alpha,\alpha}\ast u)(t;x,y).
\end{eqnarray*}
By Lemma \ref{Lemma1} we conclude that
\begin{equation*}
tu(t;x,y)+(g_1\ast S_{\alpha,\alpha-1})(t)x=\alpha[(g_1\ast u)(t;x,y)+A(g_2\ast S_{\alpha,\alpha-1}\ast u)(t;x,y)],
\end{equation*}
for any $x,y\in X$ and $t>0.$
\end{proof}

\section{Examples}

Let $-A$ be a non-negative and self-adjoint operator on the Hilbert space $X=L^2(\Omega)$ where
$\Omega\subset\mathbb R^N$ is a bounded and open set. If the operator $A$ has a compact resolvent, then $-A$ has a discrete
spectrum and its eigenvalues satisfy $0<\lambda_1\leq\lambda_2\leq\cdots\le\lambda_n\le\cdots$ with $\lim_{n\to\infty}\lambda_n=\infty.$

If $\phi_n$ denotes the normalized eigenfunction associated with $\lambda_n,$ then for all $v\in D(A)$ we have
\begin{align*}
-Av=\sum_{n=1}^\infty\lambda_n \langle v,\phi_n\rangle_{L^2(\Omega)} \phi_n.
\end{align*}
Now, consider the problem
\begin{equation}\label{eq-SAdjoint}
\left\{ \begin{array}{lcl}
\partial_t^\alpha u(t,x)&=&Au(t,x)\,\quad  t>0, \\
u(0,x)&=&u_0(x), \\
\end{array} \right.
\end{equation}
where $x\in \Omega$ and $u_0\in L^2(\Omega).$ Multiplying both sides of \eqref{eq-SAdjoint} by $\phi_n(x)$ and integrating over $\Omega$ we obtain that $u_n(t):=\langle u(t,\cdot),\phi_n(\cdot)\rangle_{L^2(\Omega)}$ is a solution of the system
\begin{equation}\label{eq-Aux-n}
\left\{ \begin{array}{lcl}
\partial_t^\alpha u_n(t)&=&-\lambda_nu_n(t)\,\quad  t>0, \\
u_n(0)&=&u_{0,n}, \\
\end{array} \right.
\end{equation}
where $u_{0,n}=\langle u_0(\cdot),\phi_n(\cdot)\rangle_{L^2(\Omega)},$ for all $n\in\mathbb{N}.$ The solution to \eqref{eq-Aux-n} is given by
\begin{equation*}
u_n(t)=E_{\alpha,1}(-\lambda_n t^\alpha)u_{0,n}=:S_{\alpha,1}^n(t)u_{0,n},
\end{equation*}
where $S_{\alpha,1}^n(t):=E_{\alpha,1}(-\lambda_n t^\alpha)$ is the resolvent family generated by $A_n:=-\lambda_n.$ According to notation in Theorem \ref{Th3.1}, we have
 $\varphi_t^n:\mathbb{R}\to\mathbb{R}$ and $\psi_t^n:\mathbb{R}\to \mathbb{R},$ are respectively, given by $\varphi_t^n(u_{0,n})=u_n(t)-u_{0,n}=S_{\alpha,1}^n(t)u_{0,n}-u_{0,n}$ and $\psi^n_t(u_{0,n})=S_{\alpha,1}^n(t)'u_{0,n}.$  By \eqref{EqCaputo5}, $S_{\alpha,1}^n(t)'=-\lambda_n S_{\alpha,\alpha}^n(t)=-\lambda_nt^{\alpha-1}E_{\alpha,\alpha}(-\lambda_n t^\alpha).$
By Theorem \ref{Th3.1},
\begin{equation*}
\alpha=\lim_{t\to 0^+} \frac{tu'_n(t)}{u_n(t)-u_{0,n}}=\lim_{t\to 0^+} \frac{-\lambda_nt^{\alpha}E_{\alpha,\alpha}(-\lambda_n t^\alpha)}{E_{\alpha,1}(-\lambda_n t^\alpha)u_{0,n}-u_{0,n}}.
\end{equation*}
Now, let $T>0$ be a fixed time. By Remark \ref{Rem1} we have
\begin{equation*}
\alpha=\frac{Tu_n(T)-(g_1\ast u_n)(T)}{(u_n\ast u_n)(T)-(g_1\ast u_n)(T)}.
\end{equation*}
By Lemma \ref{Lemma1} we have $(g_1\ast u_n)(T)=(g_1\ast S_{\alpha,1}^n)(T)u_{0,n}=S_{\alpha,2}^n(T)u_{0,n}=TE_{\alpha,2}(-\lambda_n T^\alpha)u_{0,n}$ and by \cite[Theorem 11.2]{Ha-Ma-Sa-11},
\begin{eqnarray*}
(u_n\ast u_n)(T)&=&(S_{\alpha,1}^n\ast S_{\alpha,1}^n)(T)\\
&=&\int_0^T E_{\alpha,1}(-\lambda_n(T-s)^\alpha)E_{\alpha,1}(-\lambda_n s^\alpha)u_{0,n}ds\\
&=&TE_{\alpha,2}^2(-\lambda_n T^\alpha)u_{0,n},
\end{eqnarray*}
where $E_{\alpha,2}^2(z):=\sum_{k=0}^\infty \frac{(k+1)z^k}{\Gamma(\alpha k+2)}$ is the generalized Mittag-Leffler function (see \cite{Ha-Ma-Sa-11}). Therefore,
\begin{equation*}
\alpha=\frac{TE_{\alpha,1}(-\lambda_n T^\alpha)u_{0,n}-TE_{\alpha,2}(-\lambda_n T^\alpha)u_{0,n}}{TE_{\alpha,2}^2(-\lambda_n T^\alpha)u_{0,n}-TE_{\alpha,2}(-\lambda_n T^\alpha)u_{0,n}}=\frac{E_{\alpha,1}(-\lambda_n T^\alpha)-E_{\alpha,2}(-\lambda_n T^\alpha)}{E_{\alpha,2}^2(-\lambda_n T^\alpha)-E_{\alpha,2}(-\lambda_n T^\alpha)}:=\alpha_n,
\end{equation*}
for every $T>0$ and any $n\in \mathbb{N}$ (that is, any eigenvalue $\lambda_n$ and $u_{0,n}$).

This means that, to find $\alpha$ in \eqref{eq-SAdjoint} we need know: $u_n(T)=E_{\alpha,1}(-\lambda_n T^\alpha)u_{0,n},$ $(g_1\ast u_n)(T)=E_{\alpha,2}(-\lambda_n T^\alpha)u_{0,n}$ and $(u_n\ast u_n)(T)=E_{\alpha,2}^2(-\lambda_n T^\alpha)u_{0,n}$ for any eigenvalue $\lambda_n$ of $A$ and any fixed time $T>0.$

The next Table shows a comparison between $\alpha_n,$ and $\alpha=0.2$ and $\alpha=0.4,$ for different choices of $T>0$ and $\lambda_n.$ Here, the Mittag-Leffler function has been approximated by its $N$-partial sums, with $N=50.$

\vspace{12pt}
\begin{table}[h]
\begin{center}
\begin{tabular}{| c | c | c | c || c | c | c | c |}
\hline
$\alpha$ & $T$ & $\lambda_n$ & $\alpha_n$ & $\alpha$  & $T$ & $\lambda_n$  &$\alpha_n$ \\ \hline
0.2 & 0.1 & 4 & 0.1999999998 & 0.4 & 0.1 & 4 & 0.3999999997  \\ \hline
0.2 & 0.1 & 9 & 0.1999999998 & 0.4 & 0.1 & 9 & 0.3999999998 \\ \hline
0.2 & 1 & 4 & 0.2000000002 & 0.4 & 1 & 4 & 0.3999999996 \\ \hline
0.2 & 1 & 9 & 0.2000000000 & 0.4 & 1 & 9 & 0.4000000001 \\ \hline
0.2 & 10 & 4 & 0.1999999999 & 0.4 & 10 & 4 & 0.4000000001  \\ \hline
0.2 & 10 & 9 & 0.2000000003 & 0.4 & 10 & 9 & 0.4000000002  \\ \hline
0.2 & 100 & 4 & 0.1999999999 & 0.4 & 100 & 4 & 0.3999999997  \\ \hline
0.2 & 100 & 9 & 0.2000000001 & 0.4 & 100 & 9 & 0.4000000000  \\ \hline
\end{tabular}
\end{center}
\caption{Order in Caputo fractional derivative for $0<\alpha<1.$}
\label{T1}
\end{table}

\vspace{12pt}

Now, if we consider Problem \eqref{eqCaputoII} for $1<\alpha<2,$ then the solution to the corresponding Problem \eqref{eqCaputoII} for its eigenvalue
is given by
\begin{equation*}
u_n(t)=E_{\alpha,1}(-\lambda_n t^\alpha)u_{0,n}+tE_{\alpha,2}(-\lambda_n t^\alpha)u_{1,n}=S_{\alpha,1}^n(t)u_{0,n}+S_{\alpha,2}^n(t)u_{1,n}, \quad n\in\mathbb{N}.
\end{equation*}
The notation in Theorem \ref{Th4.3} gives us $\varphi_t^n(u_{0,n},u_{1,n})=S_{\alpha,1}^n(t)u_{0,n}+S_{\alpha,2}^n(t)u_{1,n}-u_{0,n}-tu_{1,n},$ and $\psi^n_t(u_{0,n},u_{1,n})=u_n''(t)=S_{\alpha,1}^n(t)''u_{0,n}+S_{\alpha,2}^n(t)''u_{1,n}.$
By Lemma \ref{Lemma1}, $S_{\alpha,2}^n(t)=(g_1\ast S_{\alpha,1}^n)(t)$ and by \eqref{EqCaputo5} and \eqref{eqCaputo5b} we have
\begin{equation*}
S_{\alpha,2}^n(t)''=S_{\alpha,1}^n(t)'=-\lambda_nS_{\alpha,\alpha}^n(t)   \quad \mbox{ and } \quad S_{\alpha,1}^n(t)''=-\lambda_n S_{\alpha,\alpha}^n(t)'=-\lambda_n S_{\alpha,\alpha-1}^n(t).
\end{equation*}
By Theorem \ref{Th4.3} we have
\begin{equation*}
  \alpha(\alpha-1)=\lim_{t\to 0^+} \frac{ -t^2\lambda_n S_{\alpha,\alpha-1}^n(t)}{S_{\alpha,1}^n(t)u_{0,n}+S_{\alpha,2}^n(t)u_{1,n}-u_{0,n}-tu_{1,n}}.
\end{equation*}
Now, let $T>0$ be a fixed time. By Theorem \ref{Th6.1} we have
\begin{equation*}
\alpha=\frac{Tu_n(T)-(g_1\ast u_n)(T)-(g_2\ast S_{\alpha,1}^n)(T)u_{1,n}}{(S_{\alpha,1}^n\ast u_n)(T)-(g_1\ast u_n)(T)}.
\end{equation*}
By Lemma \ref{Lemma1} we have $(g_1\ast u_n)(T)=(g_1\ast S_{\alpha,1}^n)(T)u_{0,n}+(g_1\ast S_{\alpha,2}^n)(T)u_{1,n}=S_{\alpha,2}^n(T)u_{0,n}+S_{\alpha,3}^n(T)u_{1,n}=TE_{\alpha,2}(-\lambda_n T^\alpha)u_{0,n}+T^2E_{\alpha,3}(-\lambda_n T^\alpha)u_{1,n}.$
Moreover, $(g_2\ast S_{\alpha,1}^n)(T)u_{1,n}=S_{\alpha,3}^n(T)u_{1,n}=T^2E_{\alpha,3}(-\lambda_n T^\alpha)u_{1,n}$ and by \cite[Theorem 11.2]{Ha-Ma-Sa-11} we have
\begin{eqnarray*}
(S_{\alpha,1}^n\ast u_n)(T)\hspace{-0.2cm}&=&\hspace{-0.2cm}\int_0^T S_{\alpha,1}^n(T-s)u_n(s)ds\\
\hspace{-0.2cm}&=&\hspace{-0.2cm}\int_0^T E_{\alpha,1}(-\lambda_n(T-s)^\alpha)E_{\alpha,1}(-\lambda_n s^\alpha)u_{0,n}ds+\int_0^T E_{\alpha,1}(-\lambda_n(T-s)^\alpha)sE_{\alpha,2}(-\lambda_n s^\alpha)u_{1,n}ds\\
\hspace{-0.2cm}&=&\hspace{-0.2cm}TE_{\alpha,2}^2(-\lambda_n T^\alpha)u_{0,n}+T^2E_{\alpha,3}^2(-\lambda_n T^\alpha)u_{1,n},
\end{eqnarray*}
where $E_{\alpha,3}^2(z):=\sum_{k=0}^\infty \frac{(k+1)z^k}{\Gamma(\alpha k+3)}.$ Therefore
\begin{eqnarray*}\label{eqCaputo6.1}
\alpha&=&\dfrac{TE_{\alpha,1}(-\lambda_n T^\alpha)u_{0,n}+T^2E_{\alpha,2}(-\lambda_n T^\alpha)u_{1,n}-TE_{\alpha,2}(-\lambda_n T^\alpha)u_{0,n}-2T^2E_{\alpha,3}(-\lambda_n T^\alpha)u_{1,n}}{TE_{\alpha,2}^2(-\lambda_n T^\alpha)u_{0,n}+T^2E_{\alpha,3}^2(-\lambda_n T^\alpha)u_{1,n}-TE_{\alpha,2}(-\lambda_n T^\alpha)u_{0,n}-T^2E_{\alpha,3}(-\lambda_n T^\alpha)u_{1,n}}\nonumber\\
&=&\dfrac{E_{\alpha,1}(-\lambda_n T^\alpha)u_{0,n}+TE_{\alpha,2}(-\lambda_n T^\alpha)u_{1,n}-E_{\alpha,2}(-\lambda_n T^\alpha)u_{0,n}-2TE_{\alpha,3}(-\lambda_n T^\alpha)u_{1,n}}{E_{\alpha,2}^2(-\lambda_n T^\alpha)u_{0,n}+TE_{\alpha,3}^2(-\lambda_n T^\alpha)u_{1,n}-E_{\alpha,2}(-\lambda_n T^\alpha)u_{0,n}-TE_{\alpha,3}(-\lambda_n T^\alpha)u_{1,n}}\nonumber\\
&=:&\alpha_n.
\end{eqnarray*}

In the next Table we compare $\alpha_n$ and the order $\alpha=1.4$ and $\alpha=1.8$ for different choices of $T>0$ and $\lambda_n.$ For simplicity, we take $u_{0,n}=1, u_{1,n}=2.$ Here, the Mittag-Leffler function has been approximated by its $N$-partial sums, with $N=100.$

\vspace{12pt}
\begin{table}[h]
\begin{center}
\begin{tabular}{| c | c | c | c || c | c | c | c |}
\hline
$\alpha$ & $T$ & $\lambda_n$ & $\alpha_n$ & $\alpha$  & $T$ & $\lambda_n$  &$\alpha_n$ \\ \hline
1.4 & 0.5 & 1 & 1.4000000001  & 1.8 & 0.5 & 1 & 1.8000000002  \\ \hline
1.4 & 0.5 & 4 & 1.3999999999 & 1.8 & 0.5 & 4 & 1.7999999997 \\ \hline
1.4 & 1 & 1 & 1.4000000001 & 1.8 & 1 & 1 &  1.8000000002  \\ \hline
1.4 & 1 & 4 & 1.3999999999  & 1.8 & 1 & 4 & 1.8000000000 \\ \hline
1.4 & 5 & 1 & 1.4000000008 & 1.8 & 5 & 1 & 1.8000000011  \\ \hline
1.4 & 5 & 4 & 1.399959885 & 1.8 & 5 & 4 & 1.799986643  \\ \hline
\end{tabular}
\end{center}
\caption{Order in Caputo fractional derivatives for $1<\alpha<2.$}
\label{T2}
\end{table}

\vspace{12pt}

Now, we consider the fractional differential equations for the Riemann-Liouville fractional derivative \eqref{eqRLI} and \eqref{eqRLII}. Let $T>0$ be a fixed time.
We first consider $0<\alpha<1.$ The solution to the corresponding Problem \eqref{eqRLI} for its eigenvalue is given by
\begin{equation*}
u_n(t)=t^{\alpha-1}E_{\alpha,\alpha}(-\lambda_n t^\alpha)u_{0,n}=S_{\alpha,\alpha}^n(t)u_{0,n}, \quad n\in\mathbb{N}.
\end{equation*}
By Theorem \ref{Th5.4}, we have
\begin{equation*}
\alpha=\frac{Tu_n(T)}{-\lambda_n(g_1\ast S_{\alpha,\alpha}^n\ast u_n)(T)+(g_1\ast u_n)(T)}.
\end{equation*}
By Lemma \ref{Lemma1}, $(g_1\ast u_n)(T)=(g_1\ast S_{\alpha,\alpha}^n)(T)u_{0,n}=S_{\alpha,\alpha+1}^n(T)u_{0,n}=T^{\alpha}E_{\alpha,\alpha+1}(-\lambda_n T^\alpha)u_{0,n}.$
By Lemma \ref{Lemma1} and \cite[Theorem 11.2]{Ha-Ma-Sa-11} we get
\begin{eqnarray*}
(g_1\ast S_{\alpha,\alpha}^n\ast u_n)(T)&=&(S_{\alpha,\alpha+1}^n\ast u_n)(T)\\
&=&\int_0^T (T-s)^{\alpha}E_{\alpha,\alpha+1}(-\lambda_n(T-s)^\alpha)s^{\alpha-1}E_{\alpha,\alpha}(-\lambda_ns^\alpha)u_{0,n}ds\\
&=&T^{2\alpha}E_{\alpha,2\alpha+1}^2(-\lambda_n T^\alpha)u_{0,n},
\end{eqnarray*}
where $E_{\alpha,2\alpha+1}^2(z):=\sum_{k=0}^\infty \frac{(k+1)z^k}{\Gamma(\alpha k+2\alpha+1)}.$ Therefore, for any $u_{0,n},$
\begin{eqnarray*}
\alpha&=&\frac{T^{\alpha}E_{\alpha,\alpha}(-\lambda_n T^\alpha)u_{0,n}}{-\lambda_nT^{2\alpha}E_{\alpha,2\alpha+1}^2(-\lambda_n T^\alpha)u_{0,n}+T^{\alpha}E_{\alpha,\alpha+1}(-\lambda_n T^\alpha)u_{0,n}}\\
&=&\frac{E_{\alpha,\alpha}(-\lambda_n T^\alpha)}{-\lambda_nT^{\alpha}E_{\alpha,2\alpha+1}^2(-\lambda_n T^\alpha)+E_{\alpha,\alpha+1}(-\lambda_n T^\alpha)}\\
&=:&\alpha_n.
\end{eqnarray*}

The next Table compares $\alpha_n$ and the order $\alpha=0.4$ and $\alpha=0.7$ for different choices of $T>0$ and $\lambda_n.$ Here, the Mittag-Leffler function has been approximated by its $N$-partial sums, with $N=1000.$

\begin{table}[h]
\begin{center}
\begin{tabular}{| c | c | c | c || c | c | c | c |}
\hline
$\alpha$ & $T$ & $\lambda_n$ & $\alpha_n$ & $\alpha$  & $T$ & $\lambda_n$  &$\alpha_n$ \\ \hline
0.4 & 0.1 & 1 & 0.3999999998  & 0.7 & 0.1 &1 & 0.6999999993  \\ \hline
0.4 & 0.1 & 4 & 0.4000000066 & 0.7 & 0.1 & 4 & 0.7000000018 \\ \hline
0.4 & 0.5 & 1 & 0.3999999994 & 0.7 & 0.5 & 1 & 0.6999999986 \\ \hline
0.4 & 0.5 & 4 & 0.3999999709  & 0.7 & 0.5 & 4 & 0.6999998401 \\ \hline
0.4 & 1 & 1 & 0.3999999994 & 0.7 & 1 & 1 & 0.7000000000  \\ \hline
0.4 & 1 & 4 & 0.3999998780 & 0.7 & 1 & 4 & 0.6999962379  \\ \hline
\end{tabular}
\end{center}
\caption{Order in Riemann-Liouville fractional derivatives for $0<\alpha<1.$}
\label{T3}
\end{table}

\vspace{12pt}

Finally, we consider $1<\alpha<2.$ The solution to the corresponding Problem \eqref{eqRLII} for its eigenvalue is
\begin{equation*}
u_n(t)=t^{\alpha-2}E_{\alpha,\alpha-1}(-\lambda_n t^\alpha)u_{0,n}+t^{\alpha-1}E_{\alpha,\alpha}(-\lambda_n t^\alpha)u_{1,n}=S_{\alpha,\alpha-1}^n(t)u_{0,n}+S_{\alpha,\alpha}^n(t)u_{1,n}, \quad n\in\mathbb{N}.
\end{equation*}
By Theorem \ref{Th6.2},
\begin{equation*}
\alpha=\frac{Tu_n(T)+(g_1\ast S_{\alpha,\alpha-1}^n)(T)u_{0,n}}{(g_1\ast u_n)(T)-\lambda_n(g_2\ast S_{\alpha,\alpha-1}^n\ast u_n)(T)}.
\end{equation*}
The Lemma \ref{Lemma1} implies that $(g_1\ast S_{\alpha,\alpha-1}^n)(T)u_{0,n}=S_{\alpha,\alpha}^n(T)u_{0,n}=T^{\alpha-1}E_{\alpha,\alpha}(-\lambda_n T^\alpha)u_{0,n}.$ Moreover,
\begin{eqnarray*}
(g_1\ast u_n)(T)&=&(g_1\ast S_{\alpha,\alpha-1}^n)(T)u_{0,n}+(g_1\ast S_{\alpha,\alpha}^n)(T)u_{1,n}\\
&=&S_{\alpha,\alpha}^n(T)u_{0,n}+S_{\alpha,\alpha+1}^n(T)u_{1,n}\\
&=&T^{\alpha-1}E_{\alpha,\alpha}(-\lambda_n T^\alpha)u_{0,n}+T^\alpha E_{\alpha,\alpha+1}(-\lambda_n T^\alpha)u_{1,n}.
\end{eqnarray*}
Finally, by Lemma \ref{Lemma1} and \cite[Theorem 11.2]{Ha-Ma-Sa-11}  we have
\begin{eqnarray*}
(g_2\ast S_{\alpha,\alpha-1}^n\ast u_n)(T)&=&(S_{\alpha,\alpha+1}^n\ast u_n)(T)\\
&=&\int_0^T S_{\alpha,\alpha+1}(T-s)S_{\alpha,\alpha-1}(s)u_{0,n}ds+\int_0^T S_{\alpha,\alpha+1}(T-s)S_{\alpha,\alpha}(s)u_{1,n}ds\\
&=&\int_0^T (T-s)^{\alpha}E_{\alpha,\alpha+1}(-\lambda_n(T-s)^\alpha)s^{\alpha-2}E_{\alpha,\alpha-1}(-\lambda_ns^\alpha)u_{0,n}ds\\
&&+\int_0^T (T-s)^{\alpha}E_{\alpha,\alpha+1}(-\lambda_n(T-s)^\alpha)s^{\alpha-1}E_{\alpha,\alpha}(-\lambda_ns^\alpha)u_{1,n}ds\\
&=&T^{2\alpha-1}E_{\alpha,2\alpha}^2(-\lambda_n T^\alpha)u_{0,n}+T^{2\alpha}E_{\alpha,2\alpha+1}^2(-\lambda_n T^\alpha)u_{1,n}.
\end{eqnarray*}
We obtain
\begin{eqnarray*}
\alpha&=&\tfrac{T^{\alpha-1}E_{\alpha,\alpha-1}(-\lambda_n T^\alpha)u_{0,n}+T^\alpha E_{\alpha,\alpha}(-\lambda_n T^\alpha)u_{1,n}+T^{\alpha-1}E_{\alpha,\alpha}(-\lambda_n T^\alpha)u_{0,n}}{T^{\alpha-1}E_{\alpha,\alpha}(-\lambda_n T^\alpha)u_{0,n}+T^\alpha E_{\alpha,\alpha+1}(-\lambda_n T^\alpha)u_{1,n}-\lambda_nT^{2\alpha-1}E_{\alpha,2\alpha}^2(-\lambda_n T^\alpha)u_{0,n}-\lambda_nT^{2\alpha}E_{\alpha,2\alpha+1}^2(-\lambda_n T^\alpha)u_{1,n}}\\
&=&\tfrac{E_{\alpha,\alpha-1}(-\lambda_n T^\alpha)u_{0,n}+TE_{\alpha,\alpha}(-\lambda_n T^\alpha)u_{1,n}+E_{\alpha,\alpha}(-\lambda_n T^\alpha)u_{0,n}}{E_{\alpha,\alpha}(-\lambda_n T^\alpha)u_{0,n}+TE_{\alpha,\alpha+1}(-\lambda_n T^\alpha)u_{1,n}-\lambda_nT^{\alpha}E_{\alpha,2\alpha}^2(-\lambda_n T^\alpha)u_{0,n}-\lambda_nT^{\alpha+1}E_{\alpha,2\alpha+1}^2(-\lambda_n T^\alpha)u_{1,n}}\\
&=:&\alpha_n.
\end{eqnarray*}
To conclude the paper, in the next Table we compare $\alpha_n$ and the order $\alpha=1.3$ and $\alpha=1.7$ for different choices of $T>0$ and $\lambda_n.$ For simplicity, we take again $u_{0,n}=1, u_{1,n}=2.$ Moreover, the Mittag-Leffler function has been approximated by its $N$-partial sums, with $N=1000.$

\begin{table}[h]
\begin{center}
\begin{tabular}{| c | c | c | c || c | c | c | c |}
\hline
$\alpha$ & $T$ & $\lambda_n$ & $\alpha_n$ & $\alpha$  & $T$ & $\lambda_n$  &$\alpha_n$ \\ \hline
1.3 & 0.1 & 1 & 1.3000000000  & 1.7 & 0.1 & 1 & 1.7000000000  \\ \hline
1.3 & 0.1 & 4 & 1.3000000001 & 1.7 & 0.1 & 4 & 1.7000000001 \\ \hline
1.3 & 0.5 & 1 & 1.2999999999 & 1.7 & 0.5 & 1 &  1.6999999999  \\ \hline
1.3 & 0.5 & 4 & 1.2999999999  & 1.7 & 0.5 & 4 & 1.7000000000 \\ \hline
1.3 & 1 & 1 & 1.3000000002 & 1.7 & 1 & 1 & 1.7000000000  \\ \hline
1.3 & 1 & 4 & 1.3000000004 & 1.7 & 1 & 4 & 1.7000000001  \\ \hline
\end{tabular}
\end{center}
\caption{Order in Riemann-Liouville fractional derivatives for $1<\alpha<2.$}
\label{T4}
\end{table}




\begin{thebibliography}{99}


\bibitem{Ab-Al-20} L. Abadias, E. \'Alvarez, {\em Fractional Cauchy problem with memory effects,} Math. Nachr. {\bf 293} (2020), no. 10, 1846-1872.

\bibitem{Ab-Mi-15} L. Abadias, P. J. Miana, {\em A Subordination Principle on Wright Functions and Regularized Resolvent Families,} J. of Function Spaces, Volume 2015, Article ID 158145, 9 pages.

\bibitem{Al-As-20} S. Alimov, R. Ashurov, {\em Inverse problem of determining an order of the Caputo time-fractional derivative for a subdiffusion equation,} J. Inverse Ill-Posed Probl. {\bf 28} (2020), no. 5, 651-658.

\bibitem{Ar-Ba-Hi-Ne-11} W. Arendt, C. Batty, M. Hieber, F. Neubrander, {\em Vector-Valued Laplace transforms and Cauchy problems.} Monogr. Math., vol. \textbf{ 96}, Birkh\"auser,
    Basel, 2011.

\bibitem{As-Um-20} R. Ashurov, S. Umarov, {\em Determination of the order of fractional derivative for subdiffusion equations,} Fract. Calc. Appl. Anal. {\bf 23} (2020), no. 6, 1647-1662.

\bibitem{Ba-01} E. Bazhlekova, {\em Fractional evolution equations in Banach spaces,} Ph.D. thesis, Eindhoven University, 2001.

\bibitem{Ca-Pl-15} P. de Carvalho-Neto,  G. Planas, {\em Mild solutions to the time fractional Navier-Stokes equations in $\mathbb{R}^N,$} J. Differential Equations {\bf 259} (2015), no. 7, 2948-2980.

\bibitem{Cu-Pa-04} E. Cuesta, C. Palencia, {\em A numerical method for an integro-differential equations with memory in Banach spaces: Qualitative properties}, SIAM J. Numer. Anal. 41, (2003) 1232-1241.

\bibitem{Ei-Ko-04} S. Eidelman, A. Kochubei, {\em Cauchy problem for fractional diffusion equations,} J. Differential Equations {\bf 199} (2004), (2), 211-255.


\bibitem{Ha-Na-Wa-Ya-13} Y. Hatano, J. Nakagawa, S. Wang, M. Yamamoto, {\em Determination of order in fractional diffusion equation,} J. Math-for-Ind., 5A (2013), 51-57.

\bibitem{Ha-Ma-Sa-11} H. Haubold, A. Mathai, R. Saxena, {\em Mittag-Leffler Functions and Their Applications,} 2011, arXiv:0909.0230.

\bibitem{He-Li-Zo-20} J. W. He, C. Lizama, Y. Zhou, {\em The Cauchy problem for discrete-time fractional evolution equations,} J. of Computational and Applied Mathematics, {\bf 370} (2020), 112683.

\bibitem{He-Me-Po-21} H. Henr\'iquez, J. G. Mesquita, J. C. Pozo, {\em Existence of solutions of the abstract Cauchy problem of fractional order,} J. Funct. Anal. {\bf 281} (2021), no. 4, Paper No. 109028, 39 pp.



\bibitem{Ji-Ka-21} B. Jin, Y. Kian, {\em Recovering multiple fractional orders in time-fractional diffusion in an unknown medium,} Proc. A. {\bf 477} (2021), no. 2253, Paper No. 20210468, 21 pp.






\bibitem{Ji-Ru-15} B. Jin, W. Rundell, {\em A tutorial on inverse problems for anomalous diffusion processes,} Inverse Problems {\bf} 31 (2015) 035003.

\bibitem{Ka-Ru-21} B. Kaltenbacher, W. Rundell, {\em On an inverse problem of nonlinear imaging with fractional damping,} Math. Comp. 91 (2021), no. 333, 245--276.

\bibitem{Ko-Lu-19} A. Kochubei, Y. Luchko, {\em Basic FC operators and their properties,} in Handbook of fractional calculus with applications. Vol. 1, 23–46, De Gruyter, Berlin, 2019.

\bibitem{Li-Ch-Li-10} M. Li, C. Chen, F. Li, {\em On fractional powers of generators of fractional resolvent families,} J. Funct. Anal. {\bf 259} (2010) 2702-2726.

\bibitem{Li-Hu-Ya-20} Z. Li, X. Huang, M. Yamamoto, {\em A stability result for the determination of order in time-fractional diffusion equations,} J. Inverse Ill-Posed Probl. {\bf 28} (2020), no. 3, 379-388.

\bibitem{Li-Zi-Ya-19} Z. Li, Y. Liu, M. Yamamoto, {\em Inverse problems of determining parameters of the fractional partial differential equations},  in Handbook of fractional calculus with applications. Vol. 2, 431–442, De Gruyter, Berlin, 2019.

\bibitem{Li-Pe-Ji-12} K. Li, J. Peng, J. Jia, {\em Cauchy problems for fractional differential equations with Riemann-Liouville fractional derivatives,} J. Funct. Anal. 263 (2012), no. 2, 476-510.

\bibitem{Li-Zh-20} Z. Li, Z. Zhang, {\em Unique determination of fractional order and source term in a fractional diffusion equation from sparse boundary data,} Inverse Problems {\bf 36} (2020), no. 11, 115013, 20 pp.

\bibitem{Li-Li-Ya-19}  Z. Li, Y. Liu, M. Yamamoto, {\em Inverse problems of determining parameters of the fractional partial differential equations,} in Handbook of fractional calculus with applications, Vol. 2, DeGruyter (2019), p. 431-442.	



\bibitem{Li-00}  C. Lizama, {\em Regularized solutions for abstract Volterra equations,} J. Math. Anal. Appl. {\bf 243}, 278-292, (2000).

\bibitem{Li-19} C. Lizama, {\em Abstract linear fractional evolution equations,}  in Handbook of fractional calculus with applications. Vol. 2, 465–497, De Gruyter, Berlin, 2019.

\bibitem{Lu-11} S. Lukashchuk, {\em Estimation of parameters in fractional subdiffusion equations by the time integral characteristics method,} Comp. and Mathematics with Appl. {\bf 62,} 3, (2011), 834-844.

\bibitem{Pon-20} R. Ponce, {\em Asymptotic behavior of mild solutions to fractional Cauchy problems in Banach spaces,} Appl. Math. Lett. {\bf 105} (2020), 106322, 9 pp.


\bibitem{Sa-Ya-11} K. Sakamoto, M. Yamamoto, {\em Initial value/boundary value problems for fractional diffusion-wave
equations and applications to some inverse problems,} J. Math. Anal. Appl. {\bf 382} (2011) 426-447.



\bibitem{Wa-Ch-Xi-12} R. Wang, D. Chen, T. Xiao, {\em Abstract fractional Cauchy problems with almost sectorial operators,} J. Diff. Equations 252 (2012), 202-235.



\bibitem{Ya-21} M. Yamamoto, {\em Uniqueness in determining fractional orders of derivatives and initial values,} Inverse Problems {\bf 37} (2021), no. 9, Paper No. 095006, 34 pp.



\bibitem{Ze-Ch-Wa-19} X. Zheng, J. Cheng, H. Wang, {\em Uniqueness of determining the variable fractional order in variable-order time-fractional diffusion equations,} Inverse Problems {\bf 35} (2019), no. 12, 125002, 11 pp.






















\end{thebibliography}
\end{document}